\documentclass[brochure,english,12pt]{bourbaki}
\usepackage[matrix,arrow]{xy}
\usepackage{amssymb,amsfonts,amsmath,footnote}
\usepackage[francais]{babel}
\date{Novembre 2010}
\bbkannee{63\`eme ann\'ee, 2010-2011}
\bbknumero{1029}
\title{Kervaire Invariant One}
\subtitle{$\,$after M.~A.~Hill, M.~J.~Hopkins, and D.~C.~Ravenel$\,$}
\author{Haynes MILLER}
\address{Mass.~Institute of Technology\\
Department of Mathematics\\
77 Massachusetts Avenue\\
Cambridge, MA 02139--4307  -- U.S.A.}
\email{hrm@math.mit.edu}

\begin{document}

\def\colim{\mathop{\mathrm{colim}}}
\def\la#1{\mathop{\longleftarrow}\limits^{#1}}
\def\ra#1{\mathop{\longrightarrow}\limits^{#1}}
\def\be{\begin{equation}}
\def\ee{\end{equation}}

\def\VV{\mathbb{V}}
\def\MMU{\mathbb{MU}}
\def\MMO{\mathbb{MO}}
\def\HH{\mathbb{H}}
\def\SS{\mathbb{S}}
\def\LL{\mathbb{L}}
\def\LLO{\mathbb{LO}}
\def\FF{\mathbb{F}}
\def\RR{\mathbb{R}}
\def\ZZ{\mathbb{Z}}
\def\CC{\mathbb{C}}
\def\NN{\mathbb{N}}
\def\EE{\mathbb{E}}
\def\KK{\mathbb{K}}
\def\KKO{\mathbb{KO}}

\def\uZZ{\underline{\mathbb{Z}}}

\def\ev{\mathrm{ev}}
\def\map{\mathrm{map}}
\def\Hom{\mathrm{Hom}}
\def\End{\mathrm{End}}
\def\Ext{\mathrm{Ext}}
\def\Sq{\mathrm{Sq}}
\def\deg{\mathrm{deg}}
\def\id{\mathrm{id}}
\def\Th{\mathrm{\mu}}
\def\ker{\mathrm{ker}\,}
\def\Ann{\mathrm{Ann}}
\def\Mon{\mathrm{Mon}}
\def\res{\mathrm{Res}}
\def\fr{\mathrm{fr}}
\def\Sk{\mathrm{Sk}}
\def\Aut{\mathrm{Aut}}
\def\ob{\mathrm{ob}\,}
\def\Fun{\mathrm{Fun}}
\def\Stab{\mathrm{Stab}}
\def\Ind{\mathrm{Ind}}
\def\Sym{\mathrm{Sym}}
\def\op{\mathrm{op}}

\def\tensor{\otimes}

\def\cS{\mathcal{S}}
\def\cA{\mathcal{A}}
\def\cM{\mathcal{M}}
\def\cI{\mathcal{I}}
\def\cT{\mathcal{T}}
\def\cC{\mathcal{C}}
\def\cF{\mathcal{F}}

\def\hS{\widehat{S}}
\def\hW{\widehat{W}}

\def\MR{\MMU_\RR}
\def\MRn{\MMU^{(n)}}

\def\ovr{\smash{\overline{r}}}
\def\orho{\smash{\overline{\rho}}}
\def\oN{\smash{\overline{N}}}

\def\ofd{\Delta}

\maketitle

\noindent{\bf INTRODUCTION}

\bigskip
Around the year 1960 the theory of surgery was developed as part of a program
to classify manifolds of dimension greater than 4. Among the questions it 
addresses is this: 
Does every framed cobordism class contain a homotopy sphere?

Recall that a {\em framing} of a closed smooth manifold is an embedding into
Euclidean space together with a trivialization of the normal bundle. 
A good example is given by the circle embedded in $\RR^3$, with a 
framing $t$ such that the framing normal vector fields have linking number 
$\pm1$ with the circle itself. A framed manifold is {\em null-bordant} if 
it is the boundary of a framed manifold-with-boundary, and two framed manifolds
are {\em cobordant} if their difference is null-bordant.
Cobordism classes of framed $n$-manifolds form an abelian group 
$\Omega^{\fr}_n$. The class of $(S^1,t)$ in $\Omega^{\fr}_1$ is written $\eta$.

Any closed manifold of the homotopy type of $S^n$
admits a framing \cite{kervaire-milnor}, and the seemingly
absurdly ambitious question arises
of whether every class contains a manifold of the homotopy type (and hence,
by Smale, of the homeomorphism type, for $n>4$) of a sphere. 
(Classes represented by some framing of the standard $n$-sphere
form a cyclic subgroup of order 
given by the resolution of the Adams conjecture in the late 1960's.) 

A result of Pontryagin from 1950 implied that the answer had to be ``No''
in general: 
$\eta^2\neq0$ in $\Omega^{\fr}_2$, while $S^2$ is null-bordant with any 
framing.

Kervaire and Milnor \cite{kervaire-milnor} showed that the
answer is ``Yes'' unless $n=4k+2$; but that there is an obstruction, the
{\em Kervaire invariant}
\[
\kappa:\Omega^{\fr}_{4k+2}\rightarrow\ZZ/2\ZZ\,,
\] 
which vanishes on a cobordism class if and only if the class contains a 
homotopy 
sphere. Kervaire's construction \cite{kervaire} of a PL 10-manifold with no
smooth structure amounted to showing that $\kappa=0$ on $\Omega^{\fr}_{10}$. 
Pontryagin had shown that $\kappa(\eta^2)\neq0$, and in fact 
the Kervaire invariant is nontrivial on the square of any element of
Hopf invariant one. In \cite{kervaire-milnor}
Kervaire and Milnor speculated that these may be the only examples
with $\kappa\neq0$.

The Pontryagin-Thom construction establishes an isomorphism between
$\Omega_n^{\fr}$ 
and the stable homotopy group $\pi^s_n(S^0)$,
and in 1969 Browder \cite{browder} gave a homotopy theoretic interpretation
of the Kervaire invariant (Theorem \ref{browder-thm} below)
which implied that $\kappa=0$ unless
$4k+2$ is of the form $2(2^j-1)$. Homotopy theoretic calculations 
\cite{may-thesis,barratt-mahowald-tangora,barratt-jones-mahowald} 
soon muddied the waters by providing examples in dimensions 30 and 62. 
Much effort in the 1970's went into understanding the role of these 
classes and attempting inductive constructions of them in all dimensions.

The focus of this report is the following result.

\begin{theo} [Hill, Hopkins, Ravenel, 2009, \cite{hill-hopkins-ravenel}]
\label{KI1}
The Kervaire invariant \\
$\kappa:\Omega^{\fr}_{4k+2}\rightarrow\ZZ/2$
is trivial unless $4k+2=2,6,14,30,62$, or (possibly) $126$.
\end{theo}

The case $4k+2=126$ remains open.

In his proof that $\kappa$ is trivial in dimension 10, 
Kervaire availed himself of the 
current state of the art in homotopy theory (mainly work of Serre).
Kervaire and Milnor relied on contemporaneous work of Adams.
Over the intervening fifty years, further developments in homotopy theory 
have been brought to bear on the Kervaire invariant problem. In 1964, 
for example, Brown and Peterson brought spin bordism into play to show that
$\kappa$ is trivial in dimensions $8k+2$, $k>0$, 
and Browder's work used the Adams 
spectral sequence. Over the past quarter century, however, essentially no 
further progress has been made on this problem till the present work. 

Hill, Hopkins, and Ravenel (hereafter HHR) marshall three major developments 
in stable homotopy theory in their attack on the Kervaire invariant problem:
\begin{itemize}
\item The chromatic perspective based on work of
Novikov and Quillen and pioneered by Landweber, Morava, Miller, Ravenel, 
Wilson, and many more recent workers; 
\item The theory
of structured ring spectra, implemented by May and many others; and 
\item Equivariant stable homotopy theory, as developed by May and
collaborators.
\end{itemize}
The specific application of equivariant
stable homotopy theory was inspired by analogy with a fourth 
development, the motivic theory initiated by 
Voevodsky and Morel, and uses as a starting point the
theory of ``Real bordism'' investigated by Landweber, Araki,
Hu and Kriz.  
In their application of these ideas, HHR require
significant extensions of the existing state of knowledge of this subject,
and their paper provides an excellent account of the relevant parts of
equivariant stable homotopy theory.

This work was supported by the National Science Foundation Grant No. 0905950.

\section{The Kervaire invariant}

\subsection{Geometry}
Any compact smooth $n$-manifold $M$ embeds into a Euclidean space, and in 
high codimension any two embeddings are isotopic. The normal bundle $\nu$
is thus well defined up to addition of a trivial bundle. A {\em framing}
of $M$ is a bundle isomorphism $t:\nu\rightarrow M\times\RR^q$.
The Pontryagin-Thom construction is the induced
contravariant map on one-point compactifications, 
$S^{n+q}=\RR^{n+q}_+\rightarrow(M\times\RR^q)_+=M_+\wedge S^q$, giving an
element of $\pi_{n+q}(M_+\wedge S^q)$.

This group becomes independent of $q$ for 
$q>n$, and is termed the $n$th {\em stable homotopy group} 
$\pi^s_n(M_+)=\lim_{q\rightarrow\infty}\pi_{n+q}(M_+\wedge S^q)$ of $M$.
The Pontryagin-Thom construction provides 
a ``stable homotopy theory fundamental class'' $[M,t]\in\pi^s_n(M_+)$.
Composing with the map collapsing $M$ to a point gives an element of
$\pi^s_n(S^0)$. 
A bordism between framed manifolds determines a homotopy between their
Pontryagin-Thom collapse maps, and this construction gives an isomorphism
from the framed bordism group to $\pi^s_n(S^0)$. 

A framing of a manifold $M$ of dimension $4k+2$ determines additional 
structure in the cohomology of $M$. A cohomology class 
$x\in H^{2k+1}(M;\FF_2)$ is represented by a
well-defined homotopy class of maps $M_+\rightarrow K(\FF_2,2k+1)$.
When we apply the stable homotopy functor we get a map 
$\pi_*^s(M_+)\rightarrow\pi_*^s(K(\FF_2,2k+1))$. A calculation 
\cite{brown} shows that
$\pi^s_{4k+2}(K(\FF_2,2k+1))=\ZZ/2$, so the class $[M,t]$ determines an element
$q_t(x)\in\ZZ/2$. 

This {\em Kervaire form} 
$q_t:H^{2k+1}(M;\FF_2)\rightarrow\ZZ/2$ turns out to be a
{\em quadratic refinement} of the intersection pairing
$x\cdot y=\langle x\cup y,[M]\rangle$---that is to say, 
\[
q_t(x+y)=q_t(x)+q_t(y)+x\cdot y\,.
\]

In the group (under direct sum) of isomorphism classes of finite dimensional
$\FF_2$-vector spaces with nondegenerate quadratic form, 
the Kervaire forms of cobordant framed manifolds are congruent modulo the 
subgroup generated by the {\em hyperbolic} quadratic space $(H,q)$ with 
$H=\langle a,b\rangle$, $q(a)=q(b)=q(0)=0$, $q(a+b)=1$. This quotient is the 
{\em quadratic Witt group} of $\FF_2$, and is of order 2. 
The element of $\ZZ/2$ corresponding to a 
quadratic space is given by the {\em Arf invariant}, namely 
the more popular of $\{0,1\}$ as a value of $q$.
The {\em Kervaire invariant} of $(M,t)$ is then the Arf invariant of the 
quadratic space  $(H^{2k+1}(M;\FF_2),q_t)$. This defines a homomorphism
\[
\kappa:\pi^s_{4k+2}(S^0)=\Omega^{\fr}_{4k+2}\rightarrow\ZZ/2\,.
\]

\subsection{Homotopy theory}
Regarding $\kappa$ as defined on $\pi^s_{4k+2}(S^0)$ invites the question:
What is a homotopy-theoretic
interpretation of the Kervaire invariant? This was answered by Browder in a 
landmark paper \cite{browder}, in terms of the Adams spectral sequence.

Discussion of these matters is streamlined by use of the 
{\em stable homotopy category} $h\cS$, first described in Boardman's thesis
(1964). 
Its objects are called {\em spectra} and are 
designed to represent cohomology theories. There are many
choices of underlying categories of spectra, but they all lead the same
homotopy category $h\cS$, which is additive, indeed triangulated, and
symmetric monoidal (with tensor product given by the ``smash product'' 
$\wedge$). It is an analog of
the derived category of a commutative ring. There is a ``stabilization''
functor $\Sigma^\infty$ from the homotopy category of pointed CW complexes 
to $h\cS$. It sends the two point space $S^0$ to the ``sphere spectrum'' 
$\Sigma^\infty S^0=\SS$ 
which serves as the unit for the smash product. 
The suspension functor $\Sigma$ 
is given by smashing with $\Sigma^\infty S^1$.
The {\em homotopy} of a spectrum $E$ is $\pi_n(E)=[\Sigma^n\SS,E]$, so that,
for a space $X$, $\pi^s_n(X)=\pi_n(\Sigma^\infty X_+)$.
Ordinary mod $2$ cohomology of a space $X$ (which we abbreviate to $H^*(X)$)
is represented by the {\em Eilenberg-Mac Lane spectrum} 
$\HH$---$H^n(X)=[\Sigma^\infty X_+,\Sigma^n\HH]$---and,
as explained by G.~Whitehead, homology is obtained as
$H_n(X)=[\Sigma^n\SS,\Sigma^\infty X_+\wedge\HH]$. The cup product is
represented by a structure map $\HH\wedge\HH\rightarrow\HH$, making
$\HH$ into a ``ring-spectrum.'' The unit map for this ring structure,
$\SS\rightarrow\HH$, represents a generator of $\pi_0(\HH)=H^0(\SS)$.
The graded endomorphism algebra of the object $\HH$ is the well-known
{\em Steenrod algebra} $\cA$ of stable  operations on mod 2 cohomology.

Evaluation of homology gives a natural transformation (generalizing the 
{\em degree})
\[
d_X:\pi_t(X)=[\Sigma^t\SS,X]\rightarrow
\Hom_{\cA}(H^*(X),H^*(\Sigma^t\SS))
\]
which is an isomorphism if $X=\HH$. This leads to the 
{\em Adams spectral sequence}
\[
E_2^{s,t}(X;\HH)=\Ext_\cA^{s,t}(H^*(X),\FF_2)=
\Ext_\cA^{s}(H^*(X),H^*(\Sigma^t\SS))\Longrightarrow\pi_{t-s}(X)^\wedge_2
\]
for $X$ a spectrum such that $\pi_n(X)=0$ for $n\ll0$ and $H_n(X)$ is finite
dimensional for all $n$. It converges to the $2$-adic completion of the
homotopy groups of $X$. In particular, 
\[
E_2^{s,t}=\Ext_\cA^{s,t}(\FF_2,\FF_2)\Longrightarrow\pi_{t-s}(\SS)^\wedge_2\,.
\]
To this day these Ext groups remain quite mysterious overall, 
but in \cite{adams} Adams
computed them for $s\leq2$. $E_2^{0,t}$ is of course just $\FF_2$ in $t=0$.
The edge homomorphism 
\[
e:\pi_{t-1}(\SS)\rightarrow E_2^{1,t}
\]
is the ``mod 2 Hopf invariant.'' 
One interpretation of this invariant is that $e(\alpha)\neq0$ if and
only if the mapping cone $S^0\cup_\alpha e^t$ supports a nonzero 
Steenrod operation of positive degree. 
$E_2^{1,*}$ is the dual of the module of indecomposables in the Steenrod
algebra, known since Adem to be given by the Steenrod squares
$\Sq^1,\Sq^2,\Sq^4,\Sq^8,\ldots$. The element dual to $\Sq^{2^j}$
is denoted by $h_j\in E_2^{1,2^j}$. There is an element of Hopf invariant
one in $\pi_{t-1}(\SS)$ if and only if there is a nonzero permanent cycle in 
$E_2^{1,t}$. Adams \cite{adams} proved that elements of Hopf invariant one
occur only in dimensions $0$, $1$, $3$, and $7$.
This occurs because $d_2h_j=h_0h_{j-1}^2$ for $j>3$ in the Adams
spectral sequence \cite{wang}, 
a class which Adams had shown to be nontrivial.

Next, there is an edge homomorphism
\[
f:\pi_{t-2}(\SS)\rightarrow E_2^{2,t}\,,\quad t-2\neq2^j-1\,.
\]
Adams found that a basis for $E_2^{2,*}$ is given by
$\{h_i^2:i\geq0\}\cup\{h_ih_j:j>i+1\geq1\}$.
The multiplicative structure of the spectral
sequence together with computations \cite{may-thesis,mahowald-tangora} 
in $E_2^{4,*}$ imply that the only survivors to $E_3^{2,*}$ are
$h_0h_2$, $h_0h_3$, $h_2h_4$, $h_2h_5$, $h_3h_6$, and the infinite families 
$h_j^2$ and $h_1h_j$. This is the context of the following key result.

\begin{theo} [Browder, 1969, \cite{browder}]
The Kervaire invariant vanishes except in dimensions of the form
$2(2^j-1)$, where it is detected by  $h_j^2$.
\label{browder-thm}
\end{theo}

There have been subsequent simplifications of Browder's argument.
Notably,
Lannes \cite{lannes} uses manifolds with boundary
to express the Kervaire functional as a certain Hopf invariant. 
Lannes has subsequently further simplified this approach, and in
joint work with the author he has given a proof employing characteristic
numbers of manifolds with corners. Later work \cite{jones-rees}
also revealed that it is not
difficult to see, without invoking the Adams spectral sequence, that these 
are the only dimensions in which the Kervaire invariant can be nontrivial.

\subsection{Implications in homotopy theory}
Thus the Kervaire invariant problem is equivalent to the question of whether
$h_j^2$ is a permanent cycle in the Adams spectral sequence. 
Since $h_j$ is permanent for $j\leq3$, 
the product structure implies that $h_j^2$ is too---and the squares
of the framed manifolds representing the Hopf invariant one classes are
manifolds of Kervaire invariant one. Computations 
\cite{may-thesis,barratt-mahowald-tangora,barratt-jones-mahowald} showed
that $h_4^2$ and $h_5^2$ are also permanent. 
Much effort was spent searching for an inductive construction: 
If $h_j^2$ survives to an element of order 2 and square zero, 
then $h_{j+1}^2$ survives \cite{barratt-jones-mahowald}.
Parallels were developed:
Cohen, Jones, and Mahowald \cite{cohen-jones-mahowald}
defined a Kervaire invariant for oriented
$(4k+2)$-manifolds immersed in $\RR^{4k+4}$, and showed that it
does take on the value 1 in dimensions $4k+2=2(2^j-1)$. 

The theorem of HHR amounts to the assertion that $h_j^2$ supports a 
differential in the Adams spectral sequence for all $j>6$. One reason homotopy 
theorists were hoping that the Kervaire classes did survive was that they
had no idea what the targets of differentials on them might be. 
The work of HHR sheds no light on this question---though 
proving that these classes die without saying how is the genius of the paper.

Incidentally, it is known that $h_0h_2$, $h_0h_3$, and
$h_2h_4$ survive to homotopy classes
and that $h_2h_5$ does not. The fates of $h_3h_6$ and of $h_6^2$ are at present
unknown. 
All members of the remaining infinite family in 
$E_2^{2,*}$, the $h_1h_j$, do survive, to elements denoted $\eta_j$:
this is a famous theorem of Mahowald \cite{mahowald-77}. 

By the Hopf invariant one theorem, 
the operation $Sq^{2^{j+1}}$ must be trivial 
on the mapping cone of any map $S^{2^{j+1}-1}\rightarrow S^0$ if $j\geq3$.
What if we consider the mapping cone of a map 
$S^{2^{j+1}-1}\rightarrow S^0\cup_2e^1$ instead? The form of the
Adams differential $d_2h_{j+1}=h_0h_j^2$ implies that $\Sq^{2^{j+1}}$
is nonzero in the mapping cone of such a map if and only if
the composite $S^{2^{j+1}-1}\rightarrow S^0\cup_2e^1\rightarrow S^1$
has Kervaire invariant one. The result of HHR raises the
question of the minimal length \cite{christensen} of a spectrum in 
which $\Sq^{2^{j+1}}$ is nonzero.

The Kervaire invariant question is also important unstably. 
The {\em Whitehead square} $w_{n}\in\pi_{2n-1}(S^{n})$ is the composite
$S^{2n-1}\rightarrow S^{n}\vee S^{n}\rightarrow S^{n}$ of the
pinch map with the attaching map of the top cell in the torus
$S^{n}\times S^{n}$. It is null when $S^{n}$ is an $H$-space,
i.e. when $n=1,3$, or 7. For $n$ odd and not $1,3$ or 7, it is
of order 2 and generates the kernel of the suspension map 
$\pi_{2n-1}(S^{n})\rightarrow\pi_{2n}(S^{n+1})$.
Mahowald \cite{mahowald-71} proved
that if $2\theta=w_{n}$ in $\pi_{2n-1}(S^{n})$ then $\theta$ suspends 
to an element of $\pi_{2n}(S^{n+1})\cong\pi_{n-1}(\SS)$ 
of order two and Kervaire invariant one. So Browder's theorem implies
that, for $n$ odd and not 1, 3, or 7, $w_{n}$ is divisible by 2 only 
if $n+1$ is a power of 2; computations
show that it is divisible if $n=15$ or 31; and HHR prove that it is
not divisible if $n>63$.

See \cite{mahowald-71,barratt-jones-mahowald-87} 
for more about the roles elements of Kervaire 
invariant one were supposed to play in unstable homotopy theory.

\section{Outline}

In this section we state a theorem that provides the skeleton of the 
HHR proof. 

This begins with the Adams-Novikov spectral sequence, 
an analogue of the Adams spectral 
sequence in which ordinary mod $2$ homology is replaced by complex bordism,
$E_2^{*,*}(X;\MMU)\Rightarrow\pi_*(X)$. The universal Thom class is represented
by a map of spectra
$\mu:\MMU\rightarrow\HH$ which induces a map of spectral sequences
$\mu_*:E_r^{*,*}(X;\MMU)\rightarrow E_r^{*,*}(X;\HH)$.

The central innovation in the proof is the use of 
$G$-equivariant stable homotopy theory for a finite
group $G$ (e.g. the cyclic group $C_n$ of order $n$), as described in 
\S\S\ref{sec-equivariant}--\ref{sec-periodicity} below.
A $G$-spectrum $X$ has equivariant homotopy groups $\pi^G_*(X)$,
but also defines a ``homotopy fixed point spectrum'' $X^{hG}$ with its 
own homotopy groups, and there is a natural map 
$\pi^G_*(X)\rightarrow\pi_*(X^{hG})$.

\begin{theo} [HHR] 
\label{skeleton}
There exists a $C_8$-spectrum $\LL$ with the following properties.\\
{\rm (a) (Detection property)} $\pi_n(\LL)=0$ for $n$ odd, and
there are maps of spectral sequences 
\[
\xymatrix{
E^{*,*}_2(\SS;\HH) \ar@{=>}[d] & 
E^{*,*}_2(\SS;\MMU)  \ar@{=>}[d] \ar[r]^{i_*} \ar[l]_{\mu_*} &
H^*(C_8;\pi_*(\LL)) \ar@{=>}[d] \\
\pi_*(\SS) & \pi_*(\SS) \ar[l]_{=} \ar[r] & \pi_*(\LL^{hC_8})
}
\]
such that $\ker(i_*)\subseteq\ker(\mu_*)$ in $E_2^{2,2^{j}}$ for $j>6$;\\
{\rm (b) (Gap property)} $\pi_{-2}^{C_8}(\LL)=0$; \\
{\rm (c) (Fixed point property)} $\pi^{C_8}_*(\LL)\rightarrow\pi_*(\LL^{hC_8})$
is an isomorphism; and\\
{\rm (d) (Periodicity property)} 
$\pi_*(\LL^{hC_8})\cong\pi_{*+256}(\LL^{hC_8})$.
\end{theo}

\noindent
{\em Proof of Theorem \ref{KI1}}.
The facts that $\pi_{2^j-2}(\LL^{hC_8})=0$ 
(by (b), (c), and (d)) and $H^0(C_8;\pi_{2^j-1}(\LL))=0$ (by (a))
imply that any permanent cycle in $E_2^{2,2^j}(\SS;\MMU)$ maps to zero
in $H^2(C_8;\pi_{2^j}(\LL))$, and hence (by (a))
reduces to zero in $E_2^{2,2^j}(\SS;\HH)$. But since
$E_2^{1,2^j-1}(\SS;\MMU)=0$ and $E_2^{0,2^j-2}(\SS;\MMU)=0$,
any class in $\pi_{2^j-2}(\SS)$ represented by an element of
$E_2^{2,2^j}(\SS;\HH)$ must be the image of
a permanent cycle in $E_2^{2,2^j˜}(\SS;\MMU)$.
$\Box$

\medskip

The line of attack taken by HHR is suggested by an approach used in 1978 by
Ravenel \cite{ravenel} to prove an analogous result for primes $p>3$. 
We take this up in \S\ref{sec-chromatic}, as it introduces the ``chromatic''
homotopy theory underlying the current work. 
In \S\ref{sec-detection} the modifications required of Ravenel's 
method are sketched. This leads to the hope that one can take for $\LL$
the four-fold smash product $\MMU^{(4)}$ with a $C_8$ action induced
from the $C_2$ action by complex conjugation on $\MMU$. Verifying that one can
define this group action requires some technology from the modern theory of 
structured ring spectra and abstract homotopy theory, 
\S\ref{sec-ring-spectra}, and equivariant
homotopy theory, \S\ref{sec-equivariant}. 

This candidate needs to be modified if the other properties are to hold, 
however. The modification is to localize $\MMU^{(4)}$ by inverting a 
suitable equivariant homotopy class $D$. Each property puts constraints
on $D$; they are gathered in \S\ref{sec-conclusion}. 
The most severe constraint is imposed by the requirements of the 
Detection Property.

\noindent
{\em Remark.} The November, 2010, draft of HHR \cite{hill-hopkins-ravenel} 
uses $\Omega_{\mathbb{O}}$ 
for $\LL$ and $\Omega$ for $\LL^{hC_8}$. The notation used in this report 
more closely reflects the parallel between $\LL$ with its $C_8$-action
and the complex $K$-theory spectrum $\KK$ with its action of $C_2$ by complex
conjugation. The fixed point spectrum and homotopy fixed point spectrum
of that action also coincide, and are given by the orthogonal $K$-theory
spectrum  $\KKO$.
So one might write $\LL^{hC_8}=\LLO$. In any case, the spectrum $\LL$,
or $\Omega_{\mathbb{O}}$, may in time be replaced by a more tailored
construction.

\section{Chromatic homotopy theory}
\label{sec-chromatic}

\subsection{Complex cobordism and formal groups}
Much work in Algebraic Topology in the 1960s and early 1970s 
centered on developing the machinery of generalized homology theories,
satisfying the Eilenberg-Steenrod axioms with the exception of the axiom
specifying the homology of a point.
These were exemplified by $K$-theory (created by Atiyah and
Hirzebruch as
an interpretation of Bott periodicity) and cobordism theory
(Atiyah's interpretation of the Pontryagin-Thom theory of cobordism).
Novikov \cite{novikov} and Quillen \cite{quillen} 
stressed the centrality and convenience of complex 
cobordism, whose coefficient ring had been shown by Milnor
and Novikov to be a polynomial
algebra: $MU^*=MU^*(*)=\ZZ[x_1,x_2,\ldots]$, $|x_i|=-2i$. They called attention
to the class of ``complex oriented''  multiplicative cohomology theories 
$E^*(-)$: those such that the cohomology of $\CC P^\infty$ is a
power series ring over the coefficients on a single 2-dimensional 
generator---an Euler class for complex line bundles. Many of the standard
calculations of ordinary cohomology carry over to this much wider setting.
See \cite{adams-book} and \cite{ravenel-book} for this material.

A signal complication is that 
the Euler class of a tensor product is no longer the sum of the Euler
classes of the factors, but is rather given by a power series
\[
e(\lambda_1\tensor\lambda_2)=F(e(\lambda_1),e(\lambda_2))\,.
\]
Standard properties of the tensor product imply that this power series
obeys identities making it a {\em formal group}:
\[
F(x,y)=x+y+\cdots\quad,\quad
F(x,F(y,z))=F(F(x,y),z)\quad,\quad
F(x,y)=F(y,x)\,.
\]
More precisely it is a graded formal group (always of degee $2$) 
over the graded ring $E^*=E^*(*)=E_{-*}(*)=E_{-*}=\pi_{-*}(\EE)$.
The formal group associated with complex $K$-theory, for example,
is the {\em multiplicative} formal group $G_m(x,y)=x+y-uxy$ where
$u\in K^{-2}(*)$ is the Bott class, a unit by Bott periodicity. 
The spectrum $\MMU$ representing complex cobordism 
is initial among complex oriented ring spectra, and Quillen observed
that its formal group enjoys the corresponding universal property
(Lazard): there is a bijection from graded ring-homomorphisms 
$MU_*\rightarrow R_*$ to the set of graded formal groups over $R_*$
given by applying the homomorphism to the coefficients in the formal
group defined by $\MMU$. 

The integral homology $H_*(\MMU;\ZZ)$ 
also admits a modular interpretation: graded ring homomorphisms from it
are in bijection with formal groups $F(x,y)$ 
equipped with a strict isomorphism $\log(t)$
(a ``logarithm'') to the additive group $G_a(x,y)=x+y$:
$\log(x)\equiv x\!\mod\deg\, 2$ and $\log G(x,y)=\log(x)+\log(y)$.
Thus $H_*(\MMU;\ZZ)=\ZZ[m_1,m_2,\ldots]$ where $\log(t)=\sum_im_{i-1}t^i$,
and the Hurewicz map $\pi_*(\MMU)\rightarrow H_*(\MMU;\ZZ)$ classifies
the formal group $\log^{-1}(\log(x)+\log(y))$. 

Many things work just as well for $\MMU$ as for $\HH$.
For instance, we have the Adams-Novikov spectral sequence,
\[
E_2^{s,t}(X;\MMU)\Longrightarrow\pi_{t-s}(X)\,.
\]
The $E_2$ term is computable by homological algebra over the Hopf algebra 
$S_*$ that co-represents the functor sending
a commutative ring to the group (under composition) of formal powers series
$f(t)=t+b_1t^2+b_2t^3+\cdots$. As an algebra, $S_*=\ZZ[b_1,b_2,\ldots]$.
Landweber and Novikov observed that $MU_*(X)$ is 
naturally a comodule for this graded Hopf algebra (in which we declare
$|b_n|=2n$). For example, the coaction $\psi:MU_*\rightarrow S_*\tensor MU_*$
(with $X=*$) is the map co-representing the action of power series on the
set of formal groups by conjugation, ${}^fF(x,y)=f(F(f^{-1}(x),f^{-1}(y)))$. 
Then
\[
E_2^{s,t}(X;\MMU)=\Ext^{s,t}_{S_*}(\ZZ,MU_*(X))\,.
\]

The universal mod $2$ Thom class for complex vector bundles is represented by
a map $\Th:\MMU\rightarrow\HH$, and this map induces a map of spectral 
sequences
\[
\mu_*:E_r^{s,t}(X;\MMU)\rightarrow E_r^{s,t}(X;\HH)\,.
\]
The mod 2 Hopf invariant factors as 
\[
\pi_{t-1}(\SS)\rightarrow E_2^{1,t}(\SS;\MMU)\ra{\mu_*} E_2^{1,t}(\SS,\HH)\,.
\]
Novikov \cite{novikov} observed (see also \cite{miller-ravenel-wilson})
that for $j>3$ the class $h_j$ is not in the image of $\Th$. 
This solves the Hopf invariant one problem, but
does not  tell us the differential in the Adams spectral sequence 
that is nonzero on $h_j$ (though Novikov \cite{novikov} (see also
\cite{miller}) went on to use complex bordism to recover the Adams
differential as well).

\subsection{Odd primary ``Kervaire'' classes}
One may hope to address the Kervaire invariant one problem by an analogous
strategy. Unfortunately, $h_j^2$ is always in the image of $\mu_*$.
In 1978 Ravenel \cite{ravenel} turned this defect to an advantage
to prove that a certain analogue of $h_j^2$ in the classical mod $p$ 
Adams spectral sequence (for $p$ a prime larger than 3) does not survive.
The HHR gambit is a variation on Ravenel's, which we therefore describe.

The elements in question are denoted $b_j$ and they lie in
$\Ext_\cA^{2,2(p-1)p^{j+1}}(\FF_p,\FF_p)$. When $p=2$, $b_j=h_{j+1}^2$. 
The first one, $b_0$, survives to an element $\beta_1\in\pi_{2(p-1)p-2}(\SS)$.
Toda proved in 1967 that for $p$ odd, $d_{2p-1}b_1=uh_0b_0^p\neq0$ 
in the Adams spectral sequence, where $u\in\FF_p^\times$
and $h_0\in\Ext_\cA^{1,2p-2}(\FF_p,\FF_p)$ represents 
$\alpha_1\in\pi_{2p-3}(\SS)$, the first element of order $p$ in $\pi_*(\SS)$.

\begin{theo} [Ravenel, 1978, \cite{ravenel}]
For $p\geq5$ and $j\geq1$, the element $b_j$ dies in the Adams spectral 
sequence.
\label{odd-primary}
\end{theo}

The elements $h_0$ and $b_j$ are images of similarly defined elements 
in $E_2^{2,*}(\SS;\MMU)$, 
and Toda's differential forces an analogous differential there
(destroying Novikov's conjecture
\cite{novikov} that for odd primes it collapsed at $E_2$). 
Ravenel used Toda's calculation to ground an induction proving that
in the Adams-Novikov spectral sequence
\[
d_{2p-1}b_{j+1}\equiv uh_0b_j^p \mod \Ann\,b_0^{a_j}\quad,\quad
a_j=p\,\frac{p^j-1}{p-1}\quad,\quad u\in\ZZ_{(p)}^\times\,.
\]
Because of the spacing of nonzero groups at $E_2^{*,*}$,
this is the first possibly nontrivial differential, so $b_{j+1}$ 
will die as long as $h_0b_j^p\neq0$ in $E_2^{2p+1,*}$. This is 
well beyond the computable range in general. However, Ravenel proved the 
following theorem, guaranteeing that the class $b_{j+1}$ dies in the
Adams-Novikov spectral sequence:

\begin{prop} [Ravenel, 1978, \cite{ravenel}] 
For any sequence $i_0,i_1,\ldots,i_k\geq0$, the element
$h_0b_0^{i_0}b_1^{i_1}\cdots b_k^{i_k}\in E_2^{*,*}(\SS;\MMU)$ is nonzero.
\label{non-vanishing}
\end{prop}

Corresponding elements occur in $E_2^{*,*}(\SS;\HH)$, but no such theorem
holds: this is one of the advantages of $\MMU$ relative to  $\HH$. 

These classes are detected by a map
$(f,\varphi):(MU_*,S_*)\rightarrow(R_*,B_*)$ of rings acted on by Hopf 
algebras.
The Hopf algebra $B_*$ arises from a finite group $G$ as the ring $R_*^G$
of functions from $G$ to $R_*$ with Hopf algebra structure given by
$\Delta(f)(g_1,g_2)=f(g_1g_2)$, under the identification
$R_*^G\tensor_{R_*}R_*^G\ra{\cong}R_*^{G\times G}$ sending
$f_1\tensor f_2$ to $(g_1,g_2)\mapsto f_1(g_1)f_2(g_2)$. 
A graded ring homomorphism 
$f:MU_*\rightarrow R_*$ classifies a graded formal group $F$ over $R_*$. 
The map $\varphi$ arises from an action of $G$ on $F$ as a group
of strict automorphisms: formal power series $g(t)\equiv t\mod\deg\,2$ that
conjugate $F$ to itself. 
The map $\varphi:S_*\rightarrow R_*^G$ is determined by its composites
with the evaluation maps $\ev_g:R_*^G\rightarrow R_*$, and 
$\ev_g\circ\varphi$ is required to classify the formal power series $g(t)$.
With these definitions, 
\[
\Ext_{R_*^G}^{s,t}(R_*,R_*)\cong H^s(G;R_t)\,.
\]

\subsection{Automorphisms of formal groups}
One must now find a graded formal group admitting a finite group of
automorphisms with rich enough cohomology. It suffices to study graded formal
groups $F$ over graded rings of the form $R_*=A[u^{\pm1}]$, $|u|=-2$.
Such a formal group must be of the form $F(x,y)=u^{-1}F_0(ux,uy)$ for a 
formal group $F_0$ over $A$. 

Formal groups over fields are well understood. Over an algebraically
closed field of characteristic $p$, the ``height'' of a formal group 
determines it up to
isomorphism. The height is defined in terms of the self-map given by
$[p](x)$, where $[k](x)$ is defined inductively by $[0](x)=0$ and
$[k](x)=F([k-1](x),x)$. A key elementary fact is that either $[p](t)=0$
(as is the case for the additive group $G_a(x,y)=x+y$), 
or $[p](t)=f(t^{p^n})$ with 
$f'(0)\neq0$. The {\em height} is the integer $n$ (declared to be infinite
if $[p](t)=0$). 
There is a formal group of any height $n$ defined over $\FF_p$ with the
property that all endomorphisms over $\overline\FF_p$ are already defined
over $\FF_{p^n}$.  This endomorphism ring is well-studied (see e.g.
\cite{ravenel-book}), and known to
contain elements of multiplicative order $p$ exactly when $(p-1)|n$. 

Ravenel took
the first possible case: $F_0$ is formal group of height $n=p-1$
over $A=\FF_q$, $q=p^n$, such that $\Aut_{\FF_q}(F_0)$ contains an element 
of order $p$. Let $C_p$ be the subgroup generated by any such element.
Then the map 
\[
\varphi_*:\Ext_{S_*}^{*,*}(\ZZ,MU_*)\rightarrow H^*(C_p;R_*)
=E[h]\tensor R_*[b]\quad,\quad|h|=1,|b|=2
\]
sends $h_0\mapsto u^{-n}h$ and $b_j\mapsto u^{-np^{j+1}}b$ for all $j$, 
proving Proposition \ref{non-vanishing}.

We have not yet proved Theorem \ref{odd-primary}, however, 
since, for large $j$, $b_{j+1}\in E_2^{2,*}(\SS;\MMU)$ is not the only
class mapping to $b_{j+1}\in E_2^{2,*}(\SS;\HH)$. 
The claim is that $d_{2p-1}$ takes on the same value on all of them,
modulo $\ker\varphi_*$; that is, $d_{2p-1}(\ker\mu_*)\subseteq\ker\varphi_*$.
Consequently
no class in $E_2^{2,*}(\SS;\MMU)$ 
mapping to $b_{j+1}\in E_2^{2,*}(\SS;\HH)$ survives, and so
$b_{j+1}$ does not survive either. 

To see that $d_{2p-1}(\ker\mu_*)\subseteq\ker\varphi_*$, Ravenel invoked
L.~Smith's construction of a certain 
8-cell complex $\VV(2)$, with bottom cell $i:\SS\rightarrow\VV(2)$, such that
the map $\varphi:MU_*\rightarrow R_*$ factors as 
$MU_*(\SS)\rightarrow MU_*(\VV(2))\rightarrow R_*$. Thus 
$\ker i_*\subseteq\ker\varphi_*$.
He checked that $\ker\mu_*\subseteq\ker i_*$. So if $a\in\ker\mu_*$ then
$i_*d_{2p-1}a=d_{2p-1}i_*a=0$, so $d_{2p-1}a\in\ker\varphi_*$.

Smith's complex $\VV(2)$ exists only for $p>3$, and this is where that
assumption is required. In fact, at $p=3$, $b_1$ dies but $b_2$ survives
in the Adams spectral sequence, and the fate of the others is at present
unknown. It is hoped that a modification of the HHR technique will resolve 
the case $p=3$ as well. 

The chromatic approach to homotopy theory rests on the following insight,
which we owe to Morava: The automorphism groups of formal groups 
over fields control the structure of the Adams-Novikov $E_2$-term. 
Ravenel's result was the first one 
to depend upon this insight in an essential way. 

\section{The Detection Theorem}
\label{sec-detection}

In this section, we will sketch a line of thought that leads
from the ideas of \S\ref{sec-chromatic} 
to a candidate for the spectrum $\LL$ of Theorem \ref{skeleton}. 

\subsection{Hopf algebroids and actions on rings}
The starting point is finding a finite group $G$
acting (this time nontrivially) on a graded ring $R_*$ and a map
$\overline i_*:E_2^{*.*}(\SS,\MMU)\rightarrow H^*(G;R_*)$ such that 
$\ker\,\overline i_*\subseteq\ker\,\mu_*$ for $(*,*)=(2,2^j)$ for large $j$.
The construction of $\overline i_*$ requires a slight extension
of the notion of a Hopf algebra. A {\em Hopf algebroid} over a commutative
ring $K$ is a co-groupoid object in the category of commutative 
(perhaps graded) algebras over $K$. A Hopf algebroid co-represents a functor
from commutative $K$-algebras to groupoids. For example, the functor sending
a commutative graded ring $R_*$ to the groupoid of formal groups and
strict isomorphisms over $R_*$ 
is corepresentable by the Hopf algebroid $(MU_*,MU_*\tensor S_*)$, with 
structure morphisms arising from the Hopf algebra structure of $S_*$ and 
its coaction on $MU_*$. Comodules for Hopf algebroids and the corresponding
Ext groups are defined in the evident way. A comodule for $MU_*\tensor S_*$
is the same thing as an $MU_*$-module over $S_*$ and  
$\Ext_{MU_*\tensor S_*}^{*,*}(MU_*,M)=\Ext_{S_*}^{*,*}(\ZZ,M)$.

An action of a finite group on a graded ring determines
a Hopf algebroid with object ring $R_*$, morphism ring $R_*^G$,
and structure maps given as follows. The maps 
$\eta_L,\eta_R:R_*\rightarrow R_*^G$
are characterized by $\ev_g\circ\eta_L=\id$, $\ev_g\circ\eta_R=g$,
$\epsilon=\ev_1:R_*^G\rightarrow R_*$, and the diagonal map is
defined by $\Delta(f)(g_1,g_2)=f(g_1g_2)$ under the identification
$\alpha:R_*^G\tensor_{R_*}R_*^G\rightarrow R_*^{G\times G}$
such that 
\[
\ev_{g_1,g_2}\alpha(f_1\tensor f_2)=f_1(g_1)\cdot g_1f(g_2)\,.
\]
A subtlety here is that (as usual in the Hopf algebroid setup)
$R_*$ must act via $\eta_R$ on the
left factor and via $\eta_L$ on the right factor in the tensor product. 
Again, $\Ext_{R_*^G}(R_*,R_*)=H^*(G;R_*)$. 

A map $(MU_*,MU_*\tensor S_*)\rightarrow(R_*,R_*^G)$ of Hopf algebroids
classifies a graded formal group $F$ over $R_*$ along with a strict 
isomorphism $\theta_g:F\rightarrow gF$, for every $g\in G$, such that
\be
\theta_1=\id \quad\mathrm{and}\quad
\theta_{g_1g_2}=g_1\theta_{g_2}\circ\theta_{g_1}\,.
\label{mumu-to-rg}\
\ee

Construction of the relevant action begins with the involution on the 
ring spectrum $\MMU$ arising from complex conjugation. 
In homology the action is given by $\overline m_j=(-1)^jm_j$; so 
as graded commutative algebras with $C_2$ action
\be
H_*(\MMU;\ZZ)=\Sym\left(\bigoplus_{j>0}\sigma^{\tensor j}\right)
\label{muh}
\ee
where $\sigma$ is the sign representation on $\ZZ$ and $\Sym$ denotes
the symmetric algebra functor. 
In homotopy this involution classifies the formal group 
\[
\overline F(x,y)=-F(-x,-y)\,.
\]
Note that $-[-1]_F(t)$ is a strict isomorphism $F\rightarrow\overline F$. 
If $R_*$ is a graded ring with involution $x\mapsto\gamma x$, then
the map $\pi_*(\MMU)\rightarrow R_*$ classifying a formal group $F$ 
over $R_*$ is equivariant if and only if $\overline F=\gamma F$. 

\subsection{The norm}
The action of $C_2$ on $\MMU$ may be promoted to an action of a 
finite group containing $C_2$ 
via the {\em norm} construction. This construction has its origins 
in group cohomology \cite{evens} and is applicable in
any symmetric monoidal category.

For any category $\cM$, write $\cM^G$ for the category of $G$-objects in 
$\cM$, and equivariant maps. 
Let $H$ be a subgroup of $G$ of finite index. If $\cM$ has finite
coproducts, the restriction map $\res_H^G:\cM^G\rightarrow\cM^H$ has a
left adjoint $\Ind^G_H$. The object underlying $\Ind_H^GX$ is the 
coproduct of $[G:H]$ copies of the object underlying $X$. 

For example, let $\cM$ be the category of commutative $\otimes$-monoids in
a symmetric monoidal category 
$\cC$. The category $\cM$ has finite coproducts, given by the tensor product,
so $\res_H^G:\cM^G\rightarrow\cM^H$ has a left adjoint. 
A key observation of HHR is that this ``multiplicative induction'' descends
to a functor $N_H^G:\cC^H\rightarrow\cC^G$, the {\em norm} functor.

The norm functor is transitive and multiplicative:
\[
N_K^GN_H^KM\cong N_H^GM\quad,\quad 
N_H^G(M\tensor P)\cong N_H^GM\tensor N_H^GP\,.
\]
It is also distributive over coproducts. When $H$ is central in $G$ 
this distributivity admits the following description.
Let $I$ be a set and suppose given $M_i\in\cC$ for each $i\in I$.
Pick a set $F$ of orbit 
representatives for the right action of $G$ on $\map(G/H,I)$.
For $f\in F$ let $\Stab(f)=\{g\in G:fg=f\}$. Then
\be
N_H^G\left(\coprod_{i\in I}M_i\right)\cong
\coprod_{f\in F}\Ind_{\Stab(f)}^G 
N_H^{\Stab(f)}\left(\bigotimes_{g\in G/\Stab(f)}M_{f(g)}\right)\,.
\label{distributivity}
\ee

Composing the norm with the adjunction
morphism $N_H^G\res_H^GR\rightarrow R$ gives us a map
\be
\label{partial-adjoint}
\oN_H^G:\cC^H(X,\res_H^GR)\rightarrow \cC^G(N_H^GX,R)\,,
\ee
the {\em internal norm}, 
natural in the commutative $G$-monoid $R$ and $H$-object $X$.

We can apply this construction to the stable homotopy category $h\cS$.
If $H$ acts on a $\wedge$-monoid $X$ in $h\mathcal{S}$ such that 
$H_*(X;\ZZ)$ is free then 
\be
H_*(N_H^GX;\ZZ)\cong N_H^GH_*(X;\ZZ)\,.
\label{nh}
\ee

So we can form the $C_{2n}$-equivariant monoid in $h\cS$
\[
\MMU^{(n)}=N_2^{2n}\MMU\,,
\]
where $2n$ is shorthand for the cyclic group $C_{2n}$ of order $2n$. 
This is a fundamental object for HHR, where it is denoted $MU^{(C_{2n})}$. 
Neglecting the group action, 
$\MMU^{(n)}=\MMU^{\wedge n}$. Its homology is easily described
using (\ref{nh}), (\ref{muh}), and the fact that 
if $V$ is a torsion-free $\ZZ H$-module then 
$N_H^G\Sym V\cong\Sym(\Ind_H^G V)$: As a $C_{2n}$-algebra,
\be
H_*(\MMU^{(n)};\ZZ)=
\Sym\left(\bigoplus_{j>0}\Ind_2^{2n}\left(\sigma^{\tensor j}\right)\right)\,.
\label{hmu}
\ee

The homotopy of $\MMU^{(n)}$ is a graded 
ring with $C_{2n}$ action, and an equivariant map from it to a graded
$C_{2n}$-ring $R_*$ classifies a graded formal group $F$ over $R_*$ 
together with a strict isomorphism $\theta:F\rightarrow\gamma F$ 
(where $\gamma$ is a generator of $C_{2n}$) such that 
\be
\gamma^nF=\overline F \quad\mathrm{and}\quad  
\gamma^{n-1}\theta\circ\cdots\circ\gamma\theta\circ\theta=-[-1]_F\,.
\label{mu-cn}
\ee

Given such data, the formula
$\theta_{\gamma^k}=\gamma^{k-1}\theta\circ\cdots\circ\gamma\theta\circ\theta$
satisfies (\ref{mumu-to-rg}) and so defines a map
$(MU_*,MU_*\tensor S_*)\rightarrow(R_*,R_*^{C_{2n}})$. 
In the universal case, the ring is 
$MU^{(n)}_*=\pi_*(\MMU^{(n)})$; we get maps of Hopf 
algebroids
\[
\overline i:(MU_*,MU_*\tensor S_*)\rightarrow(MU^{(n)}_*,(MU_*^{(n)})^{C_{2n}})
\rightarrow(R_*,R_*^{C_{2n}})\,.
\]
In cohomology we get
\be
\overline i_*:E_2^{*,*}(\SS;\MMU)\rightarrow H^*(C_{2n};MU_*^{(n)})\rightarrow
H^*(C_{2n};R_*)\,.
\label{detection-map}
\ee

\subsection{Formal $A$-modules}
Triples $(R_*,F,\theta)$ satisfying (\ref{mu-cn})
arise from the theory of formal $A$-modules.
Let $A$ be a discrete valuation ring of characteristic zero, 
with quotient field $K$, uniformizer $\pi$, 
and residue field $A/(\pi)=\FF_q$. Over $K$, any formal group has the
form $F(x,y)=l^{-1}(l(x)+l(y))$ for a unique $l(t)\in tK[[t]]$ with 
$l'(0)=1$ (the ``logarithm'' of $F$). The formal group is defined over 
$A$ provided the Hazewinkel (\cite{hazewinkel}, 8.3) functional equation 
\[
l(t)=t+\pi^{-1}l(t^q)
\]
is satisfied. If this condition holds then for any 
$a\in A$ the power series $[a](t)=l^{-1}(al(t))$ has coefficients in $A$, 
and is an endomorphism of $F$. This construction defines a ring homomorphism
$[-]_F:A\rightarrow\End_A(F)$ which splits the natural map 
$\End_A(F)\rightarrow A$
given by $f(t)\mapsto f'(0)$: $F$ is a ``formal $A$-module.'' 
For computations, HHR use the example (\cite{hazewinkel}, 25.3.16) 
\[
l(t)=\sum_{i=0}^\infty\frac{t^{q^i}}{\pi^i}\,.
\]

Let $2n$ be a power of 2 and write $A=\ZZ_2[\zeta]$, 
where $\zeta$ is a primitive $(2n)$th root of unity. 
The element $\pi=\zeta-1$ serves as a uniformizer. The discrete valuation
ring  $A$ is
totally ramified over $\ZZ_2$, so $A/(\pi)=\FF_2$ and  $q=2$. Let
$F_0$ be a formal $A$-module, $R_*=A[u^{\pm1}]$ with $|u|=-2$,
and $F(x,y)=u^{-1}F_0(ux,uy)$. Define an action of $C_{2n}$ on
$R_*$ by letting a generator $\gamma$ act trivially on $A$ and by
$\gamma u=\zeta u$. Define 
\[
\theta(t)=\zeta^{-1}u^{-1}[\zeta](ut)\,.
\]
Then $\gamma^nF=\overline F$, 
$\theta:F\rightarrow\gamma F$ is a strict isomorphism satisfying 
(\ref{mu-cn}), and we get an equivariant ring homomorphism
$MU_*^{(n)}\rightarrow R_*$.
 
HHR take $k=3$ here, and check that for $j>6$ the map $\overline i_*$ of 
(\ref{detection-map}) (with $n=4$)
satisfies the Detection Property $\ker\overline i_*\subseteq\ker\mu_*$
for $(*,*)=(2,2^j)$. Smaller values of $k$ won't do.

We will see in \S\ref{sec-equivariant}
that the $C_{2n}$ action on the homotopy type $\MRn$
can be lifted to an action on the underlying ring-spectrum. As indicated 
in \S\ref{sec-fixed-point}
one can then form the homotopy fixed point spectrum $(\MRn)^{hC_{2n}}$,
and the middle term in (\ref{detection-map}) is the $E_2$ term of the homotopy
fixed point spectral sequence converging to $\pi_*((\MRn)^{hC_{2n}})$. 
Unfortunately $\pi_{2^j-2}((\MRn)^{hC_{2n}})$ fails to vanish 
(no matter what $n$ is).
To repair this, with $n=4$, HHR define $\LL$ to be
a certain localization of $\MMU^{(4)}$ whose homotopy is given by 
$\omega^{-1}\pi_*(\MMU^{(4)})$ for an element 
$\omega\in\pi_{2k}(\MMU^{(4)})$ which
maps to a unit in $R_*$. In \S\ref{sec-conclusion} 
we will indicate what is required of the localization
in order for such $\omega$ to exist.

The factorization (\ref{detection-map}) 
of $\overline i_*$ thus refines to a factorization
\[
E_2^{*,*}(\SS;\MMU)\ra{i_*}
H^*(C_{2n};\pi_*(\LL))\rightarrow H^*(C_{2n};R_*)\,.
\]
The fact that $\LL$ is an $\MMU$-module-spectrum implies that $i_*$ 
extends to a map of spectral sequences, and the proof of the Detection 
Property is complete. 

\section{Model categories and ring spectra} 
\label{sec-ring-spectra}

Working in a homotopy category, such as the derived category of a ring or
the stable homotopy category $h\cS$, offers a valuable simplification,
and the axiomatic framework of triangulated categories is very attractive. 
But it imposes severe restrictions on the types of construction one can make,
and at least since Quillen's fundamental document \cite{quillen-ha}  
homotopy theory has been recognized as something more than the study
of a homotopy category. Quillen wanted to say what a homotopy {\em theory}
was, and did so by means of the theory of {\em model categories}. 
Just as a homotopy type may be represented by many non-homeomorphic
spaces, so a homotopy theory may be modeled by many categorically 
non-equivalent (but ``Quillen equivalent'') model categories. A fundamental
example is given by the categories of spaces and simplicial sets. 
Quillen equivalent model categories have equivalent homotopy categories.

There are by now many models of the homotopy theory of spectra---many 
Quillen equivalent model categories of ``spectra,'' all having $h\cS$ 
as homotopy category. Some are even endowed with the structure of a 
symmetric monoidal category, in which the monoidal product descends 
to the smash product in $h\cS$. This convenience allows one to define
a ``(commutative) ring spectrum'' as a (commutative) monoid with respect
to that operation. These form model categories in their own right. 

Thom spectra, such as $\MMU$, provide basic examples in which the structure 
of a commutative ring spectrum arises explicitly as part of the construction.
When this does not happen, one can resort to an obstruction theory.
Using this approach, Goerss, Hopkins, and the author 
\cite{rezk,goerss-hopkins} 
built a commutative ring spectrum $E(F/k)$ associated to any
formal group $F$ of finite height over a perfect field $k$. The homotopy ring
is $\Lambda[u^{\pm1}]$, where $\Lambda$ denotes the Lubin-Tate ring
supporting the universal deformation of $F$. The naturality of this 
construction implies that any finite group $G$ of automorphisms of $F$
acts on the spectrum $E(F/k)$. 

Basic to homotopy theory is the construction of {\em homotopy limits},
in particular of the homotopy fixed point object of the action of a finite
group $G$. In the category of spaces, the homotopy fixed point space
of an action of $G$ on $X$ is provided by the space of equivariant maps
$X^{hG}=\map^G(EG,X)$ from a contractible CW complex $EG$ upon which $G$ acts 
freely. The map $EG\rightarrow*$ induces a map from the actual fixed point set,
$X^G\rightarrow X^{hG}$. The homotopy fixed point spectrum is the ``derived''
fixed point object, in the sense that it is the functor best approximating
$X\mapsto X^G$ that respects suitable weak equivalences. 
Similar constructions work in greater generality, and allow one
to form, for example, the homotopy fixed point ring spectrum $X^{hG}$ 
for an action of $G$ on a ring spectrum $X$.
This construction comes with a spectral sequence of the form
\[
E_2^{s,t}=H^s(G;\pi_t(X))\Longrightarrow\pi_{t-s}(X^{hG})\,.
\]

Hopkins, Mahowald, and the author explored these objects in the 
simplest nontrivial cases, when $F$ has height $p-1$ and $G$ 
contains an element of
order $p$. A motivating case occurs when $p=2$ and $F(x,y)$ is the
multiplicative formal group $G_m(x,y)=x+y-xy$
over $\FF_2$. One then finds that $E(G_m/\FF_2)^{hC_2}=KO^\wedge_2$, the
2-adic completion of orthogonal $K$-theory. Use was made in this work,
as well, of the equivariant homotopy theory of representation spheres,
to compute differentials in the homotopy fixed point spectral sequence.

Given this history it was natural for HHR to hope to define $\LL$ as
$E(F/k)$ with $F$ of height 4
and $G=C_8$ as dictated by the computations underlying the Detection Theorem. 
This hope foundered on the difficulty of the 
computations, and was replaced by the consideration of $\MMU^{(4)}$ 
and the far more elaborate 
appeal to equivariant topology required to prove properties (b), (c), and (d)
of Theorem \ref{skeleton}. 

The norm construction may be carried out in $\cS$ to produce a commutative
ring spectrum
$N_2^{2n}\MMU=\MMU^{(n)}$ with an action of $C_{2n}$, so the homotopy fixed
point set may be formed. Proving (b)--(d) for its localization 
$\LL$ requires further development of equivariant stable homotopy theory.

\section{Equivariant stable homotopy theory}
\label{sec-equivariant}

In order to construct the spectrum $\LL$ and study its properties, HHR invoke
a large collection of tools from equivariant stable homotopy theory. 
Diverse variants of this theory have been extensively developed, 
primarily by May and 
his collaborators, and HHR elect to work in the context of orthogonal 
spectra \cite{mandell-may}. At many points they require details that are
not fully articulated in existing references, and have provided in 
their paper an exhaustive account of the relevant fundamentals of this subject.

\subsection{$G$-CW complexes}
\label{subsection-gcw}

The starting point is $G$-equivariant {\em unstable} homotopy theory,
where $G$ is a finite group. 
A $G$-{\em CW complex} is a $G$-space $X$ together
with a filtration by ``skeleta'' such that $\Sk_0X$ is a $G$-set and for 
each $n\geq0$ there is a $G$-set $P_n$ and a pushout square
\[
\xymatrix{
P_n\times S^{n-1} \ar[r] \ar[d] & P_n\times D^n \ar[d] \\
\Sk_{n-1}X \ar[r] & \Sk_nX
}
\]
(where $D^n$ is the $n$-disk, $S^{n-1}$ is its boundary sphere, and 
$S^{-1}=\varnothing$), and $X=\bigcup\Sk_nX$ in the weak topology. 
For example, the unit sphere $S(V)$ of an orthogonal representation $V$
is a compact smooth $G$-manifold and hence (by a theorem of Verona) 
admits a finite $G$-CW structure. The one-point compactification $S^V$ of $V$
is the suspension of $S(V)$ and hence also admits a $G$-CW structure
(with cells increased by one in dimension from those of $S(V)$). 

When $G=C_{2n}$ and $n$ is a power of 2 we can be more explicit.  
Decompose $V$ into irreducibles---
\[
V=a\epsilon\oplus b\sigma\oplus\lambda_1\oplus\cdots\oplus\lambda_c
\]
---where $\epsilon$ denotes the trivial one-dimensional representation, 
$\sigma$ 
denotes the sign representation pulled back under the unique surjection
$C_{2n}\rightarrow C_2$, and each $\lambda_j$ is 2-dimensional. For each $j$,
$\lambda_j$ is a copy of $\CC$ with a chosen 
generator $\gamma\in C_{2n}$ acting 
by a root of unity, $\zeta_j$ (which is well defined up to complex 
conjugation). Order the $\lambda_j$'s so that if $j\leq k$ then
$\langle\zeta_j\rangle\supseteq\langle\zeta_k\rangle$ as subgroups of the
unit circle. Pick a basis $v_1,\ldots v_b$ of $b\sigma$,
and for each $j$ pick a nonzero vector $w_j\in\lambda_j$.
Define cones in $V$ as follows: 
\[
H_a=a\epsilon\quad,\quad
H_{a+i}=H_{a+i-1}\oplus\RR\langle v_i\rangle\quad
\hbox{for}\quad1\leq i\leq b\,,
\]
and for $j$ such that $1\leq j\leq c$
\[
H_{a+b+2j-1}=H_{a+b+2j-2}\oplus C_{2n}\cdot\RR_+\langle w_j\rangle\quad,\quad
H_{a+b+2j}=H_{a+b+2j-2}\oplus\lambda_j\,.
\]
Let $\Sk_0S^V=\varnothing_+=*$, and for $i>0$ let
\[
\Sk_i S^V=(H_i)_+\,.
\]
This defines a $C_{2n}$-CW filtration of $S^V$ with cells in dimensions
$0$ and $i$ with $a\leq i\leq\dim V$, each indexed by a transitive 
$C_{2n}$-set (with isotropy $C_{2n}$ in dimensions 0 and $a$,
$C_{n}$ for $a+1\leq i\leq a+b$,
and $\langle\zeta_j\rangle$ in dimensions $a+b+2j-1$ and 
$a+b+2j$).

\subsection{Orthogonal spectra}

While we will suppress discussion of technicalities regarding spectra,
it may be useful to spell out the definition of the variant adopted by HHR.

Let $G$ be a finite group. Orthogonal $G$-spectra are defined in terms 
of a certain category $\cI_G$ enriched in the category $\cT^G$ of
pointed $G$-spaces and equivariant maps. The objects of $\cI_G$ are the
finite-dimensional orthogonal representations of $G$. For $V,W\in\cI_G$,
let $O(V,W)$ be the Stiefel manifold of linear isometric embeddings from
$V$ into $W$. Let $N$ be the vector bundle over $O(V,W)$ whose fiber at
$i:V\hookrightarrow W$ is the orthogonal complement of the image of $i$ 
in $W$. Then $\cI_G(V,W)$ is defined to be the Thom space of this bundle.
The group $G$ acts on $O(V,W)$ by conjugation and $N$ is an equivariant 
vector bundle, so $\cI_G(V,W)$ receives a $G$-action. 

Let $\cT_G$ be the category of pointed
$G$-spaces and all continuous pointed maps; this is enriched over $\cT^G$.
An {\em (orthogonal) $G$-spectrum} is an enriched functor 
$X:\cI_G\rightarrow\cT_G$. It assigns to every $V\in\cI_G$ a $G$-space 
$X_V$, and to every pair $V,W\in\cI_G$ a continuous equivariant map
$\cI_G(V,W)\rightarrow\map(X_V,X_W)$. Thus for any $i\in O(V,W)$ we
receive a map $N(i)_+\wedge X_V\rightarrow X_W$; they vary continuously and
form an equivariant family. They are bonding maps for a spectrum. 
Morphisms in the enriched category $\cS_G$ of $G$-spectra are 
simply the spaces of natural transformations with $G$ acting by conjugation.
Write $\cS^G$ for
the category enriched over $\cT$ with $\cS^G(X,Y)=\cS_G(X,Y)^G$, and
$\cS_1=\cS$. 

For example, the {\em sphere $G$-spectrum} $\SS_G$ sends $W$ to its
one-point compactification $S^W$. 
The equivariant stabilization functor $\Sigma^\infty:\cT_G\rightarrow\cS_G$ 
is given on $K\in\cT_G$ by $W\mapsto S^W\wedge K$. 
Any $V\in\cI_G$ co-represents a $G$-spectrum, denoted by $S^{-V}$. Then
$S^0=S^{-0}=\SS_G$. By Yoneda, $\cS_G(S^{-V},X)=X_V$. This leads to a 
canonical presentation of any $G$-spectrum as
\be
X=\colim_V S^{-V}\wedge X_V\,.
\label{presentation}
\ee

Smash products are also handled gracefully in this context. The smash product
of two $G$-spectra $X$ and $Y$ is the left Kan extension of 
$(V,W)\mapsto X_V\wedge Y_W$ along the functor 
$\cI_G\times\cI_G\rightarrow\cI_G$ given by orthogonal direct sum. 
This makes $\cS^G$
a closed symmetric monoidal category, and equivariant associative and 
commutative ring spectra are defined using this monoidal structure. 

Let $H$ be a subgroup of $G$. There is a restriction functor
$\res_H^G:\cS^G\rightarrow\cS^H$, which has a left adjoint
\[
\Ind_H^G:\cS^H\rightarrow\cS^G\,.
\]
The ``Wirthm\"uller isomorphism'' asserts that
$\Ind_H^G$ is also right adjoint to $\res_H^G$. 

Given subgroups $H$ and $K$ there are
``double coset formulas''
\be
\res_K^G\Ind_H^GX\cong\bigvee_{g\in H\backslash G/K}\Ind_{L(g)}^K\res_{L(g)}^HX
\label{ind-res}
\ee
and
\be
(\Ind_H^GX)\wedge(\Ind_K^GY)\cong\bigvee_{g\in H\backslash G/K}
\Ind_{L(g)}^G\left((\res_{L(g)}^HX)\wedge(\res_{L(g)}^KY)\right)
\label{double-coset}
\ee
where $L(g)=H\cap gKg^{-1}$. Taking $H=G$ in (\ref{double-coset})
gives the ``Frobenius isomorphism''
\be
X\wedge\Ind_K^GY\cong\Ind_K^G(\res_K^GX\wedge Y)\,.
\label{frobenius}
\ee

\subsection{The norm}
\label{subsec-norm}

Let $\cI$ be the full subcategory of $\cI_G$ 
consisting of the trivial representations of $G$. The restriction functor
$\cS_G=\Fun(\cI_G,\cT_G)\rightarrow\Fun(\cI,\cT_G)=\Fun(G,\cS)$ 
participates in an adjoint equivalence of symmetric monoidal categories 
enriched over $\cT^G$: the category of $G$-spectra is equivalent to the
category of $G$-objects in spectra. Homotopy theoretic aspects are not
preserved by this equivalence, but the observation does allow one to easily
make constructions. For example, HHR use this device 
to construct a {\em norm} functor
\[
N_H^G:\cS^H\rightarrow\cS^G
\]
associated to an inclusion $H\subseteq G$ of finite groups.
This norm is compatible with the one described in \S\ref{sec-detection}
and with the adjoint
of the restriction functor from commutative $G$-ring-spectra to commutative
$H$-ring-spectra. Thus if $R$ is a commutative $G$-ring-spectrum,
there is a natural map of $G$-ring-spectra, a ``power operation''
\be
P:N_H^G\res_H^GR\rightarrow R\,,
\label{norm-adjoint}
\ee 
and if $R$ is a commutative $H$-ring-spectrum there is a natural ``inclusion''
\be
i:R\rightarrow\res_H^GN_H^GR
\label{inclusion}
\ee
of $G$-ring-spectra. For any $H$-spectrum $X$ and commutative 
$G$-ring-spectrum $R$ we have the {\em internal norm}
(\ref{partial-adjoint}) in $R$-cohomology, studied by Greenlees and May
\cite{greenlees-may}:
\be
\oN_H^G:(\res_H^GR)^*(X)\rightarrow R^*(N_H^GX)
\label{internal-norm}
\ee

The distributivity formula (\ref{distributivity}) for $H$ central in $G$
holds also for this norm, so for example
\be
N_H^G\left(\bigvee_{i\geq0}S^{iV}\right)\simeq
\bigvee_{f\in F}\Ind_{\Stab(f)}^GS^{|f|\Ind_H^{\Stab(f)}V}\quad,\quad
|f|=\sum_{g\in G/\Stab(f)}f(g)\,.
\label{norm-wedge}
\ee

This norm entitles HHR to view $\MRn$ as a commutative $C_{2n}$-ring-spectrum,
rather than just a $C_{2n}$-object in the stable homotopy category.
The action of $C_{2n}$ on the spectrum itself yields the structure of a 
$C_{2n}$-ring-spectrum on $D^{-1}\MRn$ 
for any equivariant map $D:S^{V}\rightarrow\MRn$, and permits construction
of the homotopy fixed point spectrum $(D^{-1}\MRn)^{hC_{2n}}$.

\subsection{Mackey functors and Bredon cohomology}

Given finite $G$-sets $P$ and $Q$, 
let  $M_G^+(P,Q)$ be the set of isomorphism classes of 
of finite $G$-sets $X$ equipped with equivariant maps $X\rightarrow P$,
$X\rightarrow Q$. These form the morphisms in a category, in which the 
composition is given by pull-back. Each $M_G^+(P,Q)$ is a commutative 
monoid under disjoint union. Formally adjoin inverses to get an abelian
group $M_G(P,Q)$. The composition is bi-additive, so it extends to the
group completions and we receive a pre-additive category $M_G$.

A {\em Mackey functor} is an additive functor from $M_G^{\op}$ to the category 
of abelian groups. For example, the ``Burnside Mackey functor'' is given by
$\widetilde A:P\mapsto M_G(P,*)$, so that $\widetilde A(G/K)$ is 
the Burnside group $A(K)$ of isomorphism classes of virtual finite $K$-sets. 

A $G$-module $N$ determines a Mackey functor $\underline N=\underline N_G$ with
\[
\underline N(P)=\map^G(P,N)\,.
\]
The morphism $P\la{p}X\ra{q}Q$ induces 
$\underline N(Q)\rightarrow\underline N(P)$
given by sending $f:Q\rightarrow N$ to 
\[
a\mapsto\sum_{px=a}f(qx)\,.
\]
The ``constant'' Mackey functor $\underline\ZZ$, for example, arises from
$\ZZ$ as a trivial $G$-module. 

The Mackey category $M_G$ arises
in equivariant stable homotopy theory as the full subcategory of $h\cS^G$
generated by the ``discrete'' spectra $\Sigma^\infty P_+$. 
Therefore for any equivariant 
cohomology theory $E_G^*$ representable by a $G$-spectrum the
functor $P\mapsto E_G^n(\Sigma^\infty P_+)$ defines a Mackey functor. 
Similarly, a $G$-spectrum $E$ represents a homology theory with
$E_n(X)=\pi_n^G(E\wedge X)$, which extends to a Mackey functor valued 
theory with $\underline E_n(X)(P)=E_n(\Sigma^\infty P_+\wedge X)$.
We have for example Mackey functor valued homotopy groups $\underline\pi_n$.
If we regard $G$-spectra of the form $P_+\wedge S^n$ (where $P$ is a finite
$G$-set) as equivariant spheres, then it is natural to extend the appellation
``homotopy group'' to $\underline\pi_n(X)(P)=[P_+\wedge S^n,X]^G$.

{\em Bredon cohomology} with coefficients in a Mackey functor $M$ is 
characterized by isomorphisms
$H_G^0(\Sigma_+^\infty A;M)=M(A)$ and 
$H_G^n(\Sigma_+^\infty A;M)=0$ for $n\neq0$, natural in $A\in M_G$,
and is represented by a $G$-spectrum $HM$. It may be computed from a $G$-CW
structure by means of a cellular cochain complex constructed in the usual way
\cite{lewis-may-mcclure}.
It follows for example that $H_G^*(X;\uZZ)=H^*(X/G;\ZZ)$.
There is of course also {\em Bredon homology} 
$H^G_*(X;M)=\pi^G_*(X\wedge HM)$. 

The construction of Eilenberg Mac Lane spectra is compatible with restriction
and induction, so the Wirthm\"uller and Frobenius isomorphisms imply for
$K\subseteq H$ that  
\be
H_*^G(\Ind_K^GX;\uZZ)\cong H_*^K(X;\uZZ)\,.
\label{homology-induction}
\ee

\section{Slice cells and the slice spectral sequence}

The key to proving the Gap and Periodicity Properties 
is the identification of an appropriate
equivariant analogue of the Postnikov system. We begin by recalling the
Postnikov system
in an appropriate form. See \cite{farjoun,hirschhorn} for details.

Let $\cT^c$ be the category of path-connected pointed spaces. For $n\geq0$ let
$\cT^c_{\leq n}$ be the class of spaces $A$ such that 
$\map_*(S^q,A)$ is contractible for all $q>n$. (It's equivalent to require
that $\pi_q(A)=0$ for all $q>n$.) Given any $X\in\cT^c$, there is a map
$p:X\rightarrow P^nX$ such that (1) $P^nX\in\cT^c_{\leq n}$ and 
(2) $p$ is a ``$\cT^c_{\leq n}$-equivalence'': 
$p^*:\map_*(P^nX,A)\ra{\simeq}\map_*(X,A)$ for any $A\in\cT^c_{\leq n}$. 
The map $p:X\rightarrow P^nX$ is the $n$th {\em Postnikov section} of $X$,
and is well defined up to a contractible choice.
Since $\cT^c_{\leq n-1}\subseteq\cT^c_{\leq n}$, 
there is a natural factorization
$X\rightarrow P^nX\rightarrow P^{n-1}X$, in which the second map may be
assumed to be a fibration. 
The fiber $P^n_nX$ of $P^nX\rightarrow P^{n-1}X$ is a $K(\pi_n(X),n)$. 
The inverse limit of this {\em Postnikov tower} is weakly equivalent to $X$. 

The equivariant analogue employed by HHR uses a carefully chosen set of 
$G$-spectra in place of spheres, rejecting the ``spheres'' 
$G/H_+\wedge S^n$ used in the construction of $G$-CW complexes in favor
of objects induced from a restricted class of representation spheres. 
Let $\rho=\rho_K$ denote the regular representation of $K$ on $\RR K$.

\begin{defi} A {\em slice cell} is a $G$-spectrum weakly equivalent to either
$\Ind_K^GS^{m\rho_K}$ (the ``regular'' case) or 
$S^{-1}\wedge \Ind_K^GS^{m\rho_K}$ (the ``irregular'' case)
for some subgroup $K$ of $G$ and some $m\in\ZZ$. The slice cell is
{\em isotropic} if $K\neq1$. The underlying 
homotopy type of a slice cell is a wedge of spheres of dimension 
$m|K|$ or $m|K|-1$, and this number is declared to be its {\em dimension}. 
\label{def-slice-cell}
\end{defi}

HHR use ``$\widehat S$'' to indicate a slice cell. 
The relationship between slice cells and equivariant spheres begins with the
elementary observation that a slice cell of dimension $n$ admits 
the structure of a finite $G$-CW complex with cells in dimensions $k$ with
\be
\begin{split}
\lfloor n/|G|\rfloor  \leq  k  \leq  n & \quad\hbox{if}\quad  n\geq0\\
n  \leq  k  \leq  \lfloor n/|G|\rfloor & \quad\hbox{if}\quad  n<0\,.
\label{slice-cell-range}
\end{split}
\ee

We get an important vanishing result if $G$ is a cyclic 2-group.

\begin{prop} [Cell Lemma]
\label{cell-lemma} 
For any regular isotropic $C_{2n}$-slice cell $\hS$, where $2n$ is a power
of $2$, we have $H^{C_{2n}}_j(\hS;\uZZ)=0$ for $-3\leq j\leq-1$.
\end{prop}

\noindent
{\em Proof.} Suppose $\hS=\Ind_{2k}^{2n}S^{m\rho_{2k}}$. 
If $m\geq0$ the result is clear. If $m<0$ we use (\ref{homology-induction})
to see that
$H^{C_{2n}}_j(\hS;\uZZ)=H^{C_{2k}}_j(S^{m\rho_{2k}};\uZZ)$. 
An equivariant form of Spanier-Whitehead duality then implies that 
$H^{C_{2k}}_j(S^{m\rho_{2k}};\uZZ)=H_{C_{2k}}^{-j}(S^{-m\rho_{2k}};\uZZ)$.
Now the $C_{2k}$-CW structure on representation spheres discussed in
\S\ref{subsection-gcw} readily implies that these groups are zero for 
$-3\leq j\leq-1$. $\Box$

\medskip
Slice cells restrict well: Given subgroups $H$ and $K$ of $G$,
we see from (\ref{ind-res}) that 
$\res_K^G\Ind_H^GS^{m\rho_H}$ is a wedge of slice cells. They induce well:
if $K\subseteq H$ then
\[
\Ind_H^G\Ind_K^HS^{m\rho_K}\cong\Ind_K^GS^{m\rho_K}
\quad,\quad
\Ind_H^G(S^{-1}\wedge\Ind_K^HS^{m\rho_K})\cong S^{-1}\wedge\Ind_K^GS^{m\rho_K}
\,.
\]

Slice cells do not have good multiplicative properties unless
$G$ is abelian. In that case (\ref{double-coset}) shows that 
\be
\Ind_H^GS^{a\rho_H}\wedge\Ind_K^GS^{b\rho_K}\cong
[G:H+K]\Ind_{H\cap K}^GS^{c\rho_{H\cap K}}
\label{smash-slice-cell}
\ee
where $c=a[H:H\cap K]+b[K:H\cap K]$, so a smash product of regular
slice cells splits as a wedge of regular slice cells. Even in the abelian case
smash products of irregular slice cells are not wedges of slice cells.

We now imitate the construction of the Postnikov system.
Let $\cS^G_{\leq n}$ be the class of objects $A$ of $\cS^G$ 
such that $\cS_G(\hS,A)$ is equivariantly contractible for all slice
cells $\hS$ with $\dim\hS>n$. (It's equivalent to require that
$[\hS,A]^G=0$ for all slice cells $\hS$ with $\dim\hS>n$.)
Given any $X\in\cS_G$, there is an equivariant map $p:X\rightarrow P^nX$ 
such that (1) $P^nX\in\cS^G_{\leq n}$ and 
(2) $p$ is an ``$\cS^G_{\leq n}$-equivalence'':
$p_*:\cS_G(P^nX,A)\rightarrow\cS_G(X,A)$ is an equivariant weak equivalence
for any $A\in\cS^G_{\leq n}$.
The map $p:X\rightarrow P^nX$ is the $n$th {\em slice section} of $X$. 
These functors assemble into a natural tower of fibrations, the 
{\em slice tower}, whose inverse limit is weakly equivalent to $X$.
The fiber of $P^nX\rightarrow P^{n-1}X$ is written
$P_n^nX$ and is called the {\em $n$-slice} of $X$.

The slice cells available in a given dimension vary with the
dimension, but (\ref{frobenius}) 
implies that smashing with $S^{m\rho_G}$ induces a bijection between
homotopy types of slice cells in dimension $t$ and dimension $t+m|G|$.
From this one deduces natural equivalences 
\be
S^{m\rho_G}\wedge P^tX\ra{\simeq}P^{t+m|G|}(S^{m\rho_G}\wedge X)
\label{rho-shift}
\ee
compatible with the projection maps. Also, since slice cells restrict
to wedges of slice cells, the restriction of the slice tower for $X$ 
to a subgroup $H\subseteq G$ is the slice tower of $\res_H^GX$. 
In particular the underlying tower of spectra is the Postnikov tower.

Unfortunately, slices can be quite complicated and may not be ``Eilenberg
Mac Lane objects'' in any reasonable sense. This is because 
$S^n$ (with trivial action) is usually not a slice cell. But $S^0$ and
$S^{-1}$ are (take $m=0$ and $K=1$ in Definition \ref{def-slice-cell}), 
and this implies that
$P_{-1}^{-1}X=H\underline\pi_{-1}(X)$. 
Zero-slices are also of the form $HM$.
The $n$-slice of an $n$-dimensional slice cell $\hS$ can be computed: 
$P_n^n\hS=H\underline\ZZ\wedge\hS$ if $\hS$ is regular,
and $P_n^n\hS=H\widetilde A\wedge\hS$ if $\hS$ is irregular.
In general, (\ref{slice-cell-range}) leads to the conclusion that if 
$X$ is an $n$-slice then $\underline\pi_k(X)=0$ unless
\be
\begin{split}
\lfloor n/|G|\rfloor\leq k\leq n & \quad\hbox{if}\quad  n\geq0\\
n\leq k\leq\lfloor(n+1)/|G|\rfloor & \quad\hbox{if}\quad  n<0\,.
\end{split}
\label{slice-range}
\ee

The {\em slice spectral sequence} is obtained by applying equivariant
homotopy to the slice tower. HHR index it as
\[
E_2^{s,t}=\pi_{t-s}^G(P^t_tX)\Longrightarrow\pi_{t-s}^G(X)\,.
\]
This indexing is in accord with that of the Atiyah-Hirzebruch spectral
sequence, but HHR prefer to display it as one does the Adams spectral
sequence, drawing $t-s$ horizontally and $s$ vertically. 
In this display, (\ref{slice-range}) implies
that the spectral sequence is nonzero only in two wedges, in the 
northeast and southwest quadrants, bounded by lines of slope 0 and $|G|-1$.

\section{The Slice Theorem}

\subsection{Purity}
Since slices are not generally Eilenberg Mac Lane objects, the $E_2$ term
of the slice spectral sequence is not generally expressible in terms of
homology. However, the key observation of HHR is that the slices of the
crucial $C_{2n}$-spectrum $\MRn$ are extremely well behaved.
Here is the definition. 

\begin{defi}
A $G$-spectrum $X$ is {\em pure} if for every $t$ there is a wedge $\hW_t$ of 
$t$-dimensional regular slice cells and an equivariant
weak equivalence $P_t^tX\simeq H\underline\ZZ\wedge\hW_t$.
A pure $G$-spectrum is {\em isotropic} if all these slice cells are isotropic.
\end{defi}

If $X$ is pure and $P_t^tX\simeq H\underline\ZZ\wedge\hW_t$, 
then in the slice spectral sequence
\[
E_2^{s,t}=\pi^G_{t-s}(P_t^tX)=H^G_{t-s}(\hW_t;\underline\ZZ)\,.
\]

Pure $G$-spectra have a number of attractive features, among them:

\begin{lemm}
\label{pure-invariance}
 If $X$ and $Y$ are pure and $f:X\rightarrow Y$ is such that $\res_1^{G}f$ 
is a weak equivalence then $f$ is a weak equivalence. \\
\end{lemm}
The proof is an induction on the order of $G$.

\begin{lemm}
If $X$ is an isotropic pure $C_{2n}$-spectrum, where $2n$ is a power of $2$,
then
\[
\pi^{C_{2n}}_j(X)=0 \quad\hbox{for}\quad -3\leq j\leq-1\,.
\]
\label{pure-gap}
\end{lemm}
For the proof, just feed the Cell Lemma \ref{cell-lemma} into the slice 
spectral sequence.

The slice filtration in general does not have good multiplicative properties,
but for pure spectra it does:

\begin{prop}
If $X$ and $Y$ are pure $G$-spectra, there is a pairing of 
slice spectral sequences $E_r(X)\tensor E_r(Y)\rightarrow E_r(X\wedge Y)$
compatible with the pairing in homotopy.
\label{slice-multiplicative}
\end{prop}

The mainspring of the HHR approach is contained in the following theorem.

\begin{theo} [Slice Theorem]
\label{slice-theorem}
If $2n$ is a power of $2$, then the $C_{2n}$-spectrum $\MRn$ is pure
and isotropic.
\end{theo}

Notice that (\ref{rho-shift}) implies that if $X$ is pure and isotropic
then so is any suspension $S^{l\rho_G}\wedge X$. This motivates the type
of localization used to construct $\LL$ from $\MRn$. If $R$ is a 
$G$-ring-spectrum
we may use an equivariant map $D:S^{l\rho_G}\rightarrow R$ to form a telescope
\be
R\rightarrow S^{-l\rho_G}\wedge R\rightarrow S^{-2l\rho_G}\wedge R\rightarrow
\cdots\,.
\label{telescope}
\ee
The direct limit is written $D^{-1}R$. When $G$ is a cyclic 2-group
and  $R$ is pure and isotropic, 
the Cell Lemma \ref{cell-lemma} implies that
$\pi^G_j(D^{-1}R)=0$ for $-3\leq j\leq-1$. So the Slice Theorem has the 
following corollary.

\begin{coro} [Gap Theorem]
Let $2n$ be a power of $2$ and let $D:S^{l\rho_{2n}}\rightarrow\MRn$ be any
equivariant map. Then
$\pi^{(2n)}_j(D^{-1}\MRn)=0$ for $-3\leq j\leq-1$.
\label{gap-theorem}
\end{coro}

The spectrum $\LL$ will be of the form $D^{-1}\MRn$ with $n=4$ and $D$
a certain map chosen to make the other parts of Theorem 
\ref{skeleton} work out. 

\subsection{Refinements---proving the Slice Theorem}
The proof of the Slice Therem is the heart of this work. HHR
begin by approximating $\MRn$ by a wedge of slice cells, by means of a 
``refinement.'' 

A {\em refinement of the homotopy} of a $G$-spectrum $X$ is a wedge
$\hW$ of slice cells together with an equivariant map 
$\alpha:\hW\rightarrow X$ for which there exists a map of graded
abelian groups $H_*(\hW;\ZZ)\rightarrow\pi_*(\hW)$ such that
\[
\xymatrix{
& H_*(\hW;\ZZ) \ar[dl]_{\cong} \ar[d] \ar[dr]^{\id}\\
\pi_*(X) & \ar[l]^{\alpha_*} \pi_*(\hW) \ar[r]_h & H_*(\hW;\ZZ)}
\]
commutes. Here $h$ denotes the Hurewicz map. If $X$ is a $G$-ring-spectrum,
the refinement is {\em multiplicative} if $\hW$ is equipped with a
$G$-ring-spectrum structure for which $\alpha$ is multiplicative.

For example, algebra generator of $\pi_*(\MMU)$ in dimension $2j$ can be chosen
to be restriction to the trivial group of a $C_2$-equivariant map
$S^{j\rho_2}\rightarrow\MMU$. It follows that 
\[
\alpha_1:\hW(1)=
\bigwedge_{j>0}\,\,\bigvee_{i\geq0}S^{ij\rho_2}\rightarrow\MMU
\]
(where the infinite smash product is the ``weak smash,'' the homotopy colimit
over finite sub smash products)
is a multiplicative refinement of homotopy. Note that $\hW(1)$ is 
regular and isotropic. 

HHR use the theory of formal groups to  define certain equivariant homotopy 
classes 
\be
\ovr_j^{(2n)}:S^{j\rho_2}\rightarrow\res_2^{2n}\MRn
\label{r-bar}
\ee
with a variety of properties making them central actors in this work.
Using the internal norm (\ref{norm-adjoint}) they lead to classes
\be
g_j=\oN_2^{2n}\ovr_j^{(2n)}:
S^{j\rho_{2n}}\simeq N_2^{2n}S^{j\rho_2}\ra{N\ovr_j}N_2^{2n}\res_2^{2n}\MRn
\ra{P}\MRn
\label{norm-r-bar}
\ee
which play an important role. When $j=2^k-1$ this composite is written
$\ofd_k^{(2n)}$ or $\ofd_k$.  

Using the $C_{2n}$-ring-spectrum structure of $\res_2^{2n}\MRn$, 
the classes $\ovr_j$ give a $C_{2n}$-ring-spectrum map
\[
\bigwedge_{j>0}\,\,\bigvee_{i\geq0}S^{ij\rho_2}\rightarrow\res_2^{2n}\MRn\,.
\]
Apply the norm to this map, and compose with the adjunction 
(\ref{norm-adjoint}) to get a map
\[
\alpha_n:\hW(n)=
N_2^{2n}\left(\bigwedge_{j>0}\,\,\bigvee_{i\geq0}S^{ij\rho_2}\right)
\rightarrow\MRn\,.
\]
The distributivity formula (\ref{norm-wedge}) implies
\[
N_2^{2n}\left(\bigvee_{i\geq0}S^{i\rho_2}\right)\cong
\bigvee_{f\in F}\Ind_{\Stab(f)}^GS^{|f|\rho_{\Stab(f)}}
\]
and this is a wedge of regular slice cells, isotropic since $C_2$ is normal
in $C_{2n}$. Thus by (\ref{smash-slice-cell}) and multiplicativity of the norm,
$\hW(n)$ is a wedge of regular isotropic slice cells. 

HHR prove, starting from (\ref{hmu}), 
that $\alpha_n$ is a multiplicative refinement of homotopy.

They develop a ``method of polynomial algebras'' filling out the following
argument. Write $A$ for the ring spectrum $\hW(n)$ and 
$I^t$ for the wedge of all the slice cells in $A$ of dimension at least $t$. 
We have a tower under $A$
\[
\cdots\rightarrow A/I^3\rightarrow A/I^2\rightarrow A/I^1=S^0
\]
The fiber $\hW_t=I^t/I^{t+1}$ in this tower is a wedge of slice cells
of dimension $t$ if $t$ is even, and is contractible of $t$ is odd.  
The tower maps to the slice tower of $R=\MMU^{(n)}$, and so extends to a
map of towers of $R$-module-spectra
\[
\xymatrix{
\cdots & \ar[r] & R\wedge_AA/I^3 \ar[r] \ar[d] & 
R\wedge_AA/I^2 \ar[r] \ar[d] & R\wedge_AA/I^1 \ar[d] \\
\cdots & \ar[r] & P^2R \ar[r] & P^1R \ar[r] & P^0R\,.
}
\]
A characterization of the slice tower implies that this is an equivalence
of towers. On the fibers, the maps are 
\[
R\wedge_A\hW_t=R\wedge_AI^t/I^{t+1}\ra{\simeq}P^t_tR\,.
\]
The action of $A$ on the quotient $I^t/I^{t+1}$ factors though the augmentation
$A\rightarrow S^0$, so 
\[
(R\wedge_AS^0)\wedge\hW_t\ra{\simeq}P^t_tR\,.
\]
We can now deduce the Slice Theorem \ref{slice-theorem} from the following key 
fact.

\begin{theo}[Reduction Theorem] Let $A\rightarrow R$ be the multiplicative
refinement of homotopy of $R=\MMU^{(n)}$ constructed above, and suppose that
$n$ is a power of $2$. Then $R\wedge_AS^0\simeq H\underline\ZZ$.
\label{reduction-theorem}
\end{theo}

The proof of the Reduction Theorem is a fairly elaborate induction, which
we will not attempt to summarize in this report. The case $n=1$ is due 
to Hu and Kriz \cite{hu-kriz}, and motivic analogues were known to
Hopkins and Morel.

\section{The Fixed Point Theorem}
\label{sec-fixed-point}

The Detection Property pertains to a {\em homotopy} fixed point set 
while the Gap Property deals with certain equivariant homotopy groups and
hence pertain to an actual fixed point set. The task of unraveling 
the relationship between these two constructions in various settings 
has been a common one in modern homotopy theory (Thomason, Carlsson, \ldots),
and this task plays a role here too.

Let $\cF$ be a {\em family} of subgroups of $G$, i.e.~a set of subgroups
closed under conjugation and passage to subgroups. An elementary construction
yields a $G$-CW complex $E\cF$ characterized up to equivariant homotopy 
equivalence by
\[
(E\cF)^H\simeq* \quad\hbox{if}\quad H\in\cF\quad,\quad
(E\cF)^H=\varnothing\quad\hbox{if}\quad H\not\in\cF\,.
\]
Denote the suspension $G$-spectra of these $G$-spaces with disjoint basepoint
adjoined by $E\cF_+$. We have the important cofibration sequences
\[
E\cF_+\rightarrow S^0\rightarrow\widetilde E\cF\,.
\]

Two families are of special importance: the family $\{1\}$, and the family
of proper subgroups. Write $E=E_G$ and $E'=E'_G$ for the corresponding 
$G$-spaces. Thus $E_G$ is the usual free contractible $G$-CW complex. 

For a $G$-spectrum $X$ and a subgroup $H\subseteq G$ one has the categorically 
defined fixed point sub-object  
\[
X^H=\map^G(G/H,X)=\map^H(*,\res_H^GX)
\]
The {\em homotopy fixed point} spectrum is $X^{hG}$, where $X^h=\map(E_G,X)$. 
The question arises of when the map $i:X\rightarrow X^h$ induced by 
$E_G\rightarrow*$ is a weak $G$-equivalence. 
The tool for attacking this question is the map of cofiber sequences
\[
\xymatrix{
E_+\wedge X \ar[r] \ar[d] & X \ar[r] \ar[d]^i & \widetilde E\wedge X \ar[d] \\
E_+\wedge X^h \ar[r] & X^h \ar[r] & \widetilde E\wedge X^h\,.
}
\]
The left vertical map is automatically an equivalence, since $\res_1^Gi$ is
a weak equivalence and $E$ is entirely built from free $G$-cells. It remains
to analyze the map $\widetilde E\wedge i$. 

This analysis is carried out using the 
``{\em geometric fixed point construction}''
\[
\Phi^HX=(E'_H\wedge\res_H^GX)^H\,.
\]
In many ways $\Phi^H$ behaves 
the way one expects a fixed point construction to behave. It preserves
weak $G$-equivalences. If $X$ is a $G$-spectrum and $T$ a $G$-space then
$\Phi^G(X\wedge T)\simeq\Phi^G(X)\wedge T^G$. In particular 
$\Phi^GS^W\simeq S^{W^G}$, so, for $H\subseteq G$, 
$\Phi^HS^{m\rho_G}\simeq S^{m[G:H]}$. It interacts well with the norm, 
as well (up to cofibrancy issues): for a $G$-spectrum $X$ and subgroup 
$H\subseteq G$, $\Phi^GN_H^GX\simeq\Phi^HX$. 

An easy induction implies that if $X$ is 
{\em geometrically free}---that is, $\Phi^HX\simeq*$
for all nontrivial subgroups $H\subseteq G$---then 
$\widetilde E\wedge X\simeq*$.

To obtain a condition under which the ``Tate spectrum'' 
$\widetilde E\wedge X^h$ is contractible, note that 
if $R$ is a $G$-ring-spectrum and $M$ an $R$-module-spectrum, then 
$\widetilde E\wedge R\simeq*$ implies $\widetilde E\wedge M\simeq*$. Since the 
$R$-module structure on $M$ induces one on $M^h$, we can conclude:

\begin{lemm}
\label{lemma-geo-free}
If the $G$-ring-spectrum $R$ is geometrically free, then 
$i:M\rightarrow M^h$ is a weak $G$-equivalence for any $R$-module-spectrum $M$.
\end{lemm}

HHR wish apply this to the $G$-ring-spectrum $\LL=D^{-1}\MRn$ (with $G=C_{2n}$,
$2n$ a power of $2$). This puts a constraint on 
$D:S^{m\rho_{2n}}\rightarrow\MRn$. To express it succinctly it is convenient
to regard $\pi_0(\cS^G(S^V,X))$ as a homotopy group of the $G$-spectrum $X$
in ``dimension'' $V$: $\pi^G_V(X)$. 
With care one can allow $V$ to be a virtual 
representation. If $X$ is a $G$-ring-spectrum then we obtain a graded ring,
graded by the real representation ring $RO(G)$. 

One has $\Phi^{(2)}\MMU\simeq\MMO$, the unoriented Thom spectrum; so 
$\Phi^{(2n)}\MRn\simeq\MMO$. Among the properties of the 
generators described in (\ref{r-bar}) is this: if we apply $\Phi^{(2n)}$
to the composite (\ref{norm-r-bar}) we get classes in $\pi_j(\MMO)$
which are zero for $j=2^k-1$ and indecomposable otherwise. So
$\Phi^{(2m)}\ofd_k^{(2m)}=0$.  

The adjunction morphism (\ref{inclusion}) 
provides us with a $C_{2m}$-ring-spectrum map
\be
\MMU^{(m)}\rightarrow\res^{2n}_{2m}N^{2n}_{2m}\MMU^{(m)}\cong
\res^{2n}_{2m}\MRn\,.
\label{adjunction-unit}
\ee
Using it we may ask that 
$\res^{2n}_{2m}D\in\pi^{(2m)}_{(ln/m)\rho_{2m}}(\res^{2n}_{2m}\MRn)$ 
be divisible by (the image of)
$\ofd^{(2m)}_k\in\pi^{(2m)}_{(2^k-1)\rho_{2m}}(\MMU^{(m)})$. Note 
that the adjunction implies that $\ofd^{(2m)}_k$ divides
$\res_{2m}^{2n}N_{2m}^{2n}\ofd^{(2m)}_k\in\pi^{(2n)}_{(2^k-1)\rho_{2n}}(\MRn)$
in this sense.

Lemma 
\ref{lemma-geo-free} implies the following theorem (in which $\LL=D^{-1}\MRn$).

\begin{theo} [Fixed Point Theorem]
\label{thm-fixed-point}
Let $2n$ be a power of $2$, and assume that 
$D\in\pi^{(2n)}_{l\rho_{2n}}(\MRn)$ has the following property: 
For every divisor $m$ of $n$, $\res_{2m}^{2n}D$ is divisible
by $\ofd_k^{(2m)}$ for some $k$. Then $\LL\rightarrow\LL^h$ is a weak
$C_{2n}$-equivalence, and 
$\pi_*^{C_{2n}}(\LL)\rightarrow\pi_*(\LL^{hC_{2n}})$ is an isomorphism.
\end{theo}

\section{The Periodicity Theorem}
\label{sec-periodicity}

Our candidate for $\LL$ is $\LL=D^{-1}\MRn$ for a suitable choice of $n$ 
(a power of 2) and $D\in\pi_{l\rho_{2n}}(\MRn)$.  
The Periodicity Property will follow from the construction of an element
of $\omega\in\pi^{(2n)}_{2nj}(\LL)$ (for $2n=8$ and $j=32$) that restricts 
to a unit in $\pi_{2nj}(\res^{(2n)}_1\LL)$. 
For: Such an element $\omega$ determines an equivariant map 
$\overline\omega:S^{2nj}\wedge\LL\rightarrow\LL$ which is a weak equivalence
of underlying spectra and hence induces a weak equivalence
$S^{2nj}\wedge\LL^{hC_{2n}}\rightarrow\LL^{hC_{2n}}$; this is the
Periodicity Property.
(By the Slice Theorem \ref{slice-theorem}, $\LL$ is pure, so by Lemma
\ref{pure-invariance} the map $\overline\omega$ is actually an equivariant weak
equivalence.)

Since we are inverting $D$, any divisor of $D$ becomes a unit in 
$RO(C_{2n})$-graded equivariant homotopy. The problem is to find a unit with
integral grading. HHR start with the class 
$\ofd^{(2n)}_1=\oN_2^{2n}\ovr_1^{(2n)}\in\pi^{(2n)}_{\rho_{2n}}(\MRn)$,
and arrange $D$ to be divisible by it.
They then work hard (and impose further restrictions on $D$)
to find a class $v\in\pi_{j(2n-\rho_{2n})}(D^{-1}\MRn)$,
for some $j$, that restricts to $1\in\pi_0(\res^{(2n)}_1D^{-1}\MRn)$.
Then $\omega=(\ofd_1^{(2n)})^jv\in\pi_{2nj}(\LL)$ 
will serve the purpose. The size of $j$
determines how many Kervaire invariant one classes are allowed to survive. 

The construction of the class $v$ begins with the ``orientation class''
\[
u_V\in H_d^G(S^V;\uZZ)=\pi^G_d(S^V\wedge H\uZZ)=\pi^G_{d-V}(H\uZZ)
\]
of an oriented $G$-representation $V$ of dimension $d$. This class satisfies
the identities
\be
u_{U\oplus V}=u_U\cdot u_V\quad,\quad
u_{\res^G_HV}=\res^G_Hu_V\quad,\quad
u_{\Ind^G_HW}=u_{\Ind^G_H\dim W}\cdot\oN_H^Gu_W
\label{orientation-ids}
\ee
in which $\oN_H^G$ refers to the internal norm described in
(\ref{partial-adjoint}).
In particular, $u_V\in\pi^G_{d-V}(H\uZZ)$ restricts to $1\in\pi_0(H\ZZ)$.
Note that twice any representation admits an orientation.

First consider the ``sign'' representation $\sigma=\sigma_{2n}$ of $C_{2n}$, 
the pullback of the usual sign representation under the unique surjection
$C_{2n}\rightarrow C_2$. We have the orientation class 
$u_{2\sigma}\in\pi_2^{(2n)}(S^{2\sigma}\wedge H\uZZ)$, and we 
wonder first whether it, or some power of it, lifts to a class in 
$\pi^{(2n)}_2(S^{2\sigma}\wedge\MRn)$. 

The slice spectral sequence provides potential obstructions to lifting 
a class from $H\uZZ=P^0\MRn$ to $\MRn$. Since the class of interest is not
in integral grading we have to apply $RO(G)$-graded
homotopy to the slice tower. There results family of spectral sequences
indexed by $V$, 
\[
E_2^{s,t+V}=\pi^G_{t-s+V}(P^{d+t}_{d+t}X)\Longrightarrow\pi^G_{t-s+V}(X)\,.
\]
If the spectrum $X$ is pure, 
the multiplicative structure ties them all together. 

Thus we can regard $u_{2\sigma}$ as an element of $E_2^{0,2-2\sigma}$. 
HHR need some other classes to express differentials on this class.
For any $G$-representation $V$, the inclusion of the trivial
subspace defines a $G$-map $S^0\rightarrow S^V$. If $V^G\neq0$ 
then this $G$-map is null-homotopic, but in general it is not;
it defines an ``Euler class'' 
\[
a_V\in\pi_0^G(S^V)=\pi^G_{-V}(S^0)\,.
\]

Write $a$ for the Hurewicz image of $a_\sigma$:
$a\in\pi^{(2n)}_0(S^\sigma\wedge H\uZZ)=\pi^{(2n)}_{-\sigma}(H\uZZ)$;
so $a\in E_2^{1,1-\sigma}$. It survives to
the $\MRn$ Hurewicz image of $a_\sigma$ in $\pi^{(2n)}_{-\sigma}(\MRn)$.

Finally, let  $b=1\wedge a_{\orho}:S^1\rightarrow S^{\rho_{2n}}$, 
where $\orho$ denotes
the reduced regular representation of $C_{2n}$. Composing $b^j$ with the class
$g_j$ defined in (\ref{norm-r-bar}) gives $g_jb^j\in\pi^{(2n)}_j(\MRn)$,
represented by a class $f_j\in E_2^{(2n-1)j,2nj}$
(because $S^{j\rho_{2n}}$ is a slice cell of dimension $2nj$).

No power of $u_{2\sigma}$ 
survives the slice spectral sequence, by virtue of
the following computation 
(due in case $n=1$ to Araki (unpublished) and Hu and Kriz \cite{hu-kriz}).

\begin{prop} Let $2n$ be a power of $2$.
In the $RO(C_{2n})$-graded slice spectral sequence for $\MRn$,
$u_{2\sigma}^{2^{k-1}}$ survives to $E_r$ with $r=1+2n(2^k-1)$, and 
\[
d_ru_{2\sigma}^{2^{k-1}}=a^{2^k}f_{2^k-1}\neq0\,.
\]
\end{prop}
 
These very differentials imply the survival of certain multiples 
of these classes, however.
The element $\ofd_k=\oN_2^{2n}\ovr_{2^k-1}^{(2n)}$ 
(cf (\ref{norm-r-bar})) defines a map 
$S^{(2^k-1)\rho_{2n}}\wedge\MRn\rightarrow\MRn$, via the $G$-ring-spectrum
structure of $\MRn$, which is compatible with the slice tower and gives
a map of slice spectral sequences. 
A simple identity 
shows
that the target of the differential on $u_{2\sigma}^{2^k}$ is killed by
$\ofd_k$. This is the last differential that can be nonzero on
$\ofd_k u_{2\sigma}^{2^k}$, so that class 
survives in the slice spectral sequence.
This implies that the image of $u_{2\sigma}^{2^k}$ survives in the localization
of the slice spectral sequence obtained by inverting $\ofd_k$, and so 
$u_{2\sigma}^{2^{k'}}$ does too for $k'\geq k$.


The class $\ofd_1^{(2n)}$ involves the regular representation, 
not the sign representation.
Fortunately, the internal norm (\ref{partial-adjoint})
provides a relation between their orientation 
classes. Since $\rho_2=1\oplus\sigma_2$, $u_{2\rho_2}=u_{2\sigma_2}$.
Using the equations $\Ind_H^G\rho_H=\rho_G$ and 
$\Ind_n^{2n}1=1\oplus\sigma_{2n}$, one sees that 
\[
u_{2\rho_4}=u_{2\sigma_4}^4\cdot\oN_2^4u_{2\rho_2}\quad,\quad
u_{2\rho_8}=
u_{2\sigma_8}^8\cdot\oN_4^8u_{2\rho_4}
\]
and so on. Thus 
\[
u_{2\rho_8}^{2^k}=u_{2\sigma_8}^{2^{k+3}}\cdot
\oN_4^8\bigl(u_{2\sigma_4}^{2^{k+2}}\cdot\oN_2^4 u_{2\sigma_2}^{2^k}\bigr)
\]
will survive in the localized slice spectral sequence 
provided that
\be
\label{periodicity-requirements}
\begin{split}
\ofd_{k_3}^{(8)}\quad\hbox{divides}\quad  \res_8^8D & 
\quad\hbox{for some}\quad k_3\leq k+3 \\
\ofd_{k_2}^{(4)}\quad\hbox{divides}\quad \res_4^8D & 
\quad\hbox{for some}\quad k_2\leq k+2\\
\ofd_{k_1}^{(2)}\quad\hbox{divides}\quad \res_2^8D & 
\quad\hbox{for some}\quad k_1\leq k\,.
\end{split}
\ee
Any element $v\in\pi^{(8)}_{2^{k+1}(8-\rho_8)}(D^{-1}\MMU^{(4)})$ 
represented by $u_{2\rho_8}^{2^k}$ will restrict to $1$ in the group 
$\pi_0(\res_1^{(8)}D^{-1}\MMU^{(4)})$. 
Taking \
\[
\omega=(\ofd_1^{(8)})^{2^{k+1}}v\in\pi_{2^{k+4}}(\LL)
\] 
will then give us the Periodicity Property.

\section{Wrap-up}
\label{sec-conclusion}

\noindent
{\em Proof of Theorem \ref{skeleton}}. The spectrum $\LL$ is $D^{-1}\MRn$.
We have now accumulated several requirements on $2n$ (a power of 2) and $D$:
\begin{itemize}
\item The Gap Property requires that $D\in\pi^{(2n)}_{l\rho_{2n}}(\MRn)$ for
some $l$.
\item The Fixed Point Property requires that for each divisor $m$ of $n$ there
exists $k$ such that $\ofd_k^{(2m)}$ divides $\res_{2m}^{2n}D$.
\item The Periodicity Property requires, for $n=4$, that $D$ satisfies 
(\ref{periodicity-requirements}) and that $\ofd_1^{(8)}$ divides $D$.
\item The Detection Property requires that $n\geq4$ and for $n=4$ that $D$ 
satisfies (\ref{periodicity-requirements}) for certain values
of $k_1$, $k_2$, $k_3$.
\end{itemize}

A calculation verifies that this last requirement is met by  
$k_1=4$, $k_2=2$, $k_3=1$ and no smaller values. 
Recalling that $\ofd_k^{(2m)}$ divides 
$\res_{2m}^{2n}N_{2m}^{2n}\ofd_k^{(2m)}$, we can take $k=4$ and
\[
D=N_2^8\ofd^{(2)}_4\cdot N_4^8\ofd^{(4)}_2\cdot\ofd^{(8)}_1
\,\in\,\pi_{19\rho_8}^{(8)}(\MMU^{(4)})\,.
\]
Then $v\in\pi^{(8)}_{2^{k+1}(8-\rho_{8})}(D^{-1}\MMU^{(4)})$ 
with  $2^{k+1}=32$, so $\omega=(\ofd_1^{(8)})^{32}v\in\pi_{256}(\LL)$.
$\Box$

\end{document}